\documentclass[11pt]{amsart}
\usepackage{fullpage}
\usepackage{amsmath, amssymb, amsthm, verbatim}
\usepackage[foot]{amsaddr}

\usepackage{graphicx}
\usepackage{color}
\usepackage[colorlinks=true,linkcolor=red,citecolor=blue,urlcolor=cyan]{hyperref}
\usepackage{cleveref}

\usepackage{float}

\newtheorem{thm}{Theorem}[section]
\newtheorem{prop}[thm]{Proposition}
\theoremstyle{remark}

\graphicspath{{./figs/}} 

\title{Harmonic functions on finitely-connected tori}\thanks{C.-Y. Kao acknowledges partial support from NSF grant DMS-2208373. B. Osting acknowledges partial support from NSF DMS 17-52202 and DMS 21-36198. \'E. Oudet was partially supported by the project ANR-18-CE40-0013 SHAPO financed by the French Agence Nationale de la Recherche (ANR) and by the Institut Universitaire de France.}

\author{Chiu-Yen Kao}
\address{Department of Mathematical Sciences, Claremont McKenna College, Claremont, CA }
\email{ckao@cmc.edu}

\author{Braxton Osting}
\address{Department of Mathematics, University of Utah, Salt Lake City, UT}
\email{osting@math.utah.edu}

\author{\'Edouard Oudet}
\address{LJK, Universit\'e Grenoble Alpes, France}
\email{edouard.oudet@imag.fr}

\subjclass[2020]{
30F15,  
31A25,  
35C10,  
65N25.  
}
\keywords{Harmonic function; Laplace equation; finitely-connected torus; doubly-periodic domain; elliptic function; Weierstrass elliptic function; Steklov eigenvalue}
\date{\today}

\begin{document}

\begin{abstract} 
In this paper, we prove a Logarithmic Conjugation Theorem on finitely-connected tori. The theorem states that a harmonic function can be written as the real part of a function whose derivative is analytic and a finite sum of terms involving the logarithm of the modulus of a modified Weierstrass sigma function. We implement the method using arbitrary precision and use the result to find approximate solutions to the Laplace problem and Steklov eigenvalue problem. Using a posteriori estimation, we show that the solution of the Laplace problem on a torus with a few circular holes has error less than $10^{-100}$ using a few hundred degrees of freedom and the Steklov eigenvalues have similar error. 
\end{abstract}

\maketitle
\section{Introduction}
Harmonic functions satisfying the Laplace equation, $\Delta u = 0$, 
arise in many physical applications, including 
potential flow in fluid dynamics, 
the stationary solution of heat conduction, and 
electrostatics in the absence of charges, to name just a few. 
Efficient and robust numerical approaches to solving the Laplace equation on a general domain with different boundary conditions are crucial for understanding the aforementioned applications. 
In this paper, \emph{we are particularly interested in solving the Laplace equation on finitely-connected tori}, which serves as a model problem for the study of 
heat or electrical conduction in the exterior of a periodic lattice of inclusions with prescribed temperature or 
for fluid flow through a doubly periodic array of obstacles.

\subsubsection*{Harmonic functions.}
It is well-known that every harmonic function $u$ on a simply-connected domain $\Omega \subset \mathbb C$ can be written as the real part of an analytic function, $f(z)$, 
\begin{equation}
u(z) = \Re f(z). 
\end{equation}
For finitely-connected domains, the analogous result is known as the \emph{Logarithmic Conjugation Theorem} \cite{Axler_1986,Trefethen_2018}. Let $\Omega \subset \mathbb C$ be a finitely-connected region which means that $\mathbb C \setminus \Omega$ has only finitely many bounded connected components, $\{ K_j \}_{j \in [b]}$ with $b \in \mathbb N \setminus \{0\}$. For each $j \in [b]$, let $a_j$ be a point in $K_j$. If $u$ is a harmonic function on $\Omega$, then there exists an analytic function $f$ on $\Omega$ and real numbers $c_j$, $j \in [b]$, such that 
\begin{equation}
u(z) = \Re f(z) + \sum_{j \in [b]} c_j \log | z - a_j |, 
\qquad \qquad z \in \Omega.
\end{equation}

\emph{Our main result is to extend the Logarithmic Conjugation Theorem to finitely-connected tori.} 
We consider a torus 
${\mathbb T}_\omega = \mathbb C / L_\omega$, 
where 
$L_\omega = 2 \omega_1 \mathbb Z + 2 \omega_2 \mathbb Z$ is a lattice and  
$\omega = (\omega_1,\omega_2) \in \mathbb C^2$ are half-periods, assumed not to be colinear. 
Let 
\begin{equation} \label{e:Omega}
\Omega = {\mathbb T}_\omega \setminus \cup_{j \in [b]} K_j
\end{equation}
denote the finitely-connected torus after removing $b \in \mathbb N \setminus \{0\}$ disjoint, connected compact sets $\{ K_j \}_{j \in [b]}$, with smooth boundary. We also introduce the parallelogram (fundamental domain) 
\begin{equation} \label{e:O}
\mathcal P = \left\{2 \omega_1 x + 2 \omega_2 y \in \mathbb{C} \colon (x,y) \in [0,1]^2  \right\} \setminus  \cup_{j \in [b]} K_{j}. 
\end{equation}
Note that $\Omega$ is obtained from $\mathcal P$ after identification of opposite sides. 
Recall that a meromorphic, doubly-periodic function is called an \emph{elliptic} function. Let 
\begin{equation} \label{e:sigma_hat}
\hat\sigma(z,\omega) = e^{-\frac{1}{2} \gamma_2 z^2 - \frac{1}{2}\pi |z|^2/A} \sigma(z,\omega)
\end{equation}
denote the modified Weierstrass sigma function \cite{Haldane_2018}, where 
$\sigma(z,\omega)$ is the 
Weierstrass sigma function, 
$\gamma_2 = \gamma_2 (\omega) \in \mathbb C$ is a lattice invariant, and 
$A = area(\mathbb T_\omega)$. 
We further discuss $\hat\sigma(z,\omega)$ in \cref{s:Weierstrass}, but for now just note that it is a non-holomorphic, function with a pole of order 2 at $z=0$ such that 
$|\hat\sigma(z,\omega)|$ is doubly-periodic. 

\begin{thm}
\label{t:LogConjThmTorus} 
Let $\Omega$ and $\mathcal P$ be defined as in \eqref{e:Omega} and \eqref{e:O}. 
For each $j \in [b]$, let $a_j$ be a point in $K_j$. 
If $u$ is a harmonic function on $\Omega$ (equivalently, harmonic and doubly-periodic on $\mathcal P$), then there exists an analytic function $\hat f$ on $\mathcal P$ and real numbers $c_j$, $j \in [b]$, satisfying $\sum_{j \in [b]} c_j = 0$, such that $\hat f'$ is elliptic and
\begin{equation} \label{e:TorusSeries}
u(z) = \Re \hat f(z) + \sum_{j \in [b]} c_j \log | \hat \sigma \left(z - a_j, \omega \right) |, 
\qquad \qquad z \in \Omega.
\end{equation}
If there is only one connected boundary component (i.e., $b=1$), then $c_1 = 0$ and $u(z) = \Re \hat f(z)$. 
\end{thm}

A proof of \cref{t:LogConjThmTorus} is given in \cref{s:LogConjThmProof}. 
We comment that the result in \cref{t:LogConjThmTorus} differs from the Logarithmic Conjugation Theorem for finitely-connected domains in several important ways. 
First, the modified Weierstrass sigma function, $ \log | \hat \sigma( z,\omega) |$, plays the role of $\log |z|$.
Secondly, and perhaps surprisingly, while the derivative $\hat f'$ is elliptic, the function $\hat f$ cannot always be taken to be elliptic.

\subsubsection*{Computing harmonic functions on finitely-connected tori.}
There are a variety of methods for computing harmonic functions on finitely-connected tori, including integral equation methods with multipole acceleration \cite{barnett2018unified} and the finite element method \cite{guedes1990preprocessing}.  In our approach, we are inspired by the work in \cite{Trefethen_2018} to use \cref{t:LogConjThmTorus} to represent doubly-periodic harmonic functions using a series solution.  
Let 
$\wp(z) = \wp(z,\omega)$ denote the Weierstrass elliptic function, 
$\wp^{(k)}(z,\omega)$, denote the $k$-th derivative, and 
$\hat \zeta(z) = \hat \zeta(z,\omega)$ denote the ``modified''  Weierstrass zeta function that is doubly-periodic; these will be defined in \cref{s:Weierstrass}.

\begin{thm}
\label{t:HarmonicPerDecomp}
Let $\Omega$ be a  finitely-connected torus as in \eqref{e:Omega}. For each $j \in [b]$, let $a_j$ be a point in $K_j$. If $u$ is a harmonic function on $\Omega$, then there exists  a constant $C \in \mathbb{R}$ and real coefficients  $(a_{j,k})$, $(b_{j,k})$ and $(c_j)$ such that
{\small
\begin{equation} 
\label{e:HarmonicRep}
\begin{split}
u(z) = C +  \sum_{j\in [b]} \Big[ & a_{j,-1} \Re \hat{\zeta}(z - a_j) + 
                                    b_{j,-1} \Im \hat{\zeta}(z - a_j) 
            + \sum_{k\ge0} a_{j,k}  \Re \wp^{(k)}(z-a_j) + 
             \sum_{k\ge0} b_{j,k}  \Im \wp^{(k)}(z-a_j) \\
           &  +  c_j \log | \hat \sigma \left( z - a_j \right) | \Big]
\end{split}
\end{equation}}
where $\sum_{j\in [b]} c_j = 0$.
\end{thm}

A proof of \cref{t:HarmonicPerDecomp} is given in \cref{s:LogConjThmProof}. We have chosen to represent the elliptic function $f'$ using a sum of Weierstrass functions. Similar representations have been used to find doubly periodic solutions in several applications, including 
doubly-periodic stress distributions in perforated plates \cite{linkov1999complex},  
solitary wave solutions to a nonlinear wave equation \cite{chen2006weierstrass} and nonlinear Schr{\"o}dinger equation \cite{el2016constructing}, 
lowest-Landau-level wavefunctions on the torus \cite{Haldane_2018}, 
and simulation of oil recovery \cite{astafev2014simulation}. Other representations for elliptic functions are possible, including as $R(\wp) + \wp' S(\wp)$ for some rational functions $R$ and $S$ \cite{busam2009complex}. 

\begin{figure}[t!]
\centering
\includegraphics[width=0.8\textwidth,trim = {0 3cm 0 3cm},clip]{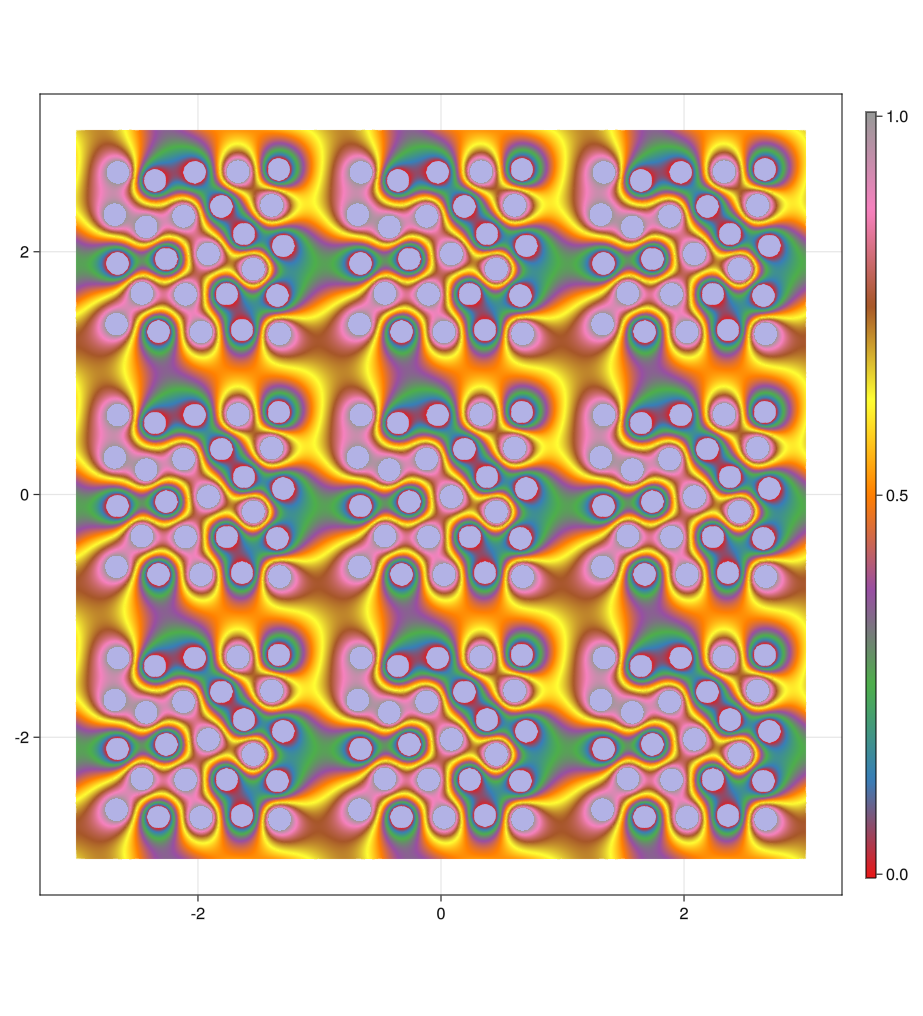}
\caption{Approximate solution to the Laplace equation on a square torus with 25 disks removed.  Dirichlet boundary conditions equal to 0 or 1 are imposed on the boundary of each disk. The computational domain is $[-1,1]^2$ and 9 copies are displayed to emphasize periodicity. 
See \cref{s:Laplace} for more details.}
\label{f:multiholes}
\end{figure}

\bigskip

In \cref{s:examples}, we use the series representation \eqref{e:HarmonicRep} to solve the Laplace problem 
\begin{subequations}
\label{e:Laplace}
\begin{align}
& \Delta u = 0  && \textrm{in \ } \Omega \\ 
& u = f && \textrm{on \ } \partial \Omega = \cup_{j\in [b]} \partial K_j,
\end{align}
\end{subequations}
where $f \in L^2(\partial \Omega)$ is given  and the Steklov eigenvalue problem
\begin{subequations}
\label{e:Steklov}
\begin{align}
& \Delta u = 0  && \textrm{in \ } \Omega \\ 
& \partial_n u = \sigma u && \textrm{on \ } \partial \Omega = \cup_{j\in [b]} \partial K_j.
\end{align}
\end{subequations}
As in  \cite{Trefethen_2018}, the series solution \eqref{e:HarmonicRep} are not convergent series. The coefficients depend on the truncation of the sum (in $k$). 
For the Laplace problem, by the maximum principle, the accuracy of the solution can be computed by looking at the error on the boundary, $\sup_{x\in \partial \Omega} |u(x) - f(x)|$. 
For the Steklov problem, we bound the error in the eigenvalues using an 
a posteriori estimate \cite{Bogosel_2016,Fox_1967}. 
We implement the proposed numerical method in Julia using arbitrary precision and use the result to find approximate solutions to the Laplace problem and Steklov eigenvalue problem. For a few circular holes, the solution of the Laplace problem has error less than $10^{-100}$ using a few hundred degrees of freedom and the Steklov eigenvalues have similar error. 
We show the solution to the Laplace problem with 25 disks removed in \cref{f:multiholes}. The spectral accuracy is also demonstrated for non-convex holes in \cref{f:LaplaceFlowers}.  

We conclude in \cref{s:disc} with a discussion. 

\section{Weierstrass elliptic functions} 
\label{s:Weierstrass}

Here we recall some background material on Weierstrass elliptic functions and establish notation used in the paper. Excellent references include  \cite{busam2009complex,Haldane_2018,pastras2020weierstrass}. 

We consider the lattice 
$$
L_\omega = 2 \omega_1 \mathbb Z + 2 \omega_2 \mathbb Z,
$$
where 
$\omega = (\omega_1,\omega_2) \in \mathbb C^2$ are half-periods, assumed not to be colinear. 
A function $f \colon \mathbb C \to \mathbb C $ is said to be \emph{doubly-periodic} if it satisfies
\begin{align*}
f(z+2\omega_{1}) & =f(z)\\
f(z+2\omega_{2}) & =f(z)
\end{align*}
for all $z\in\mathbb{C}$. 
A function is said to be \emph{elliptic} if it is meromorphic and doubly-periodic. 
An example of an elliptic function is the Weierstrass elliptic function
\[
\wp(z,\omega) := \frac{1}{z^{2}}+\sum_{\ell \in L_\omega \setminus \{0\}} \ \left(\frac{1}{\left(z-\ell \right)^{2}}-\frac{1}{\left(\ell\right)^{2}}\right).
\]
The subtraction of the last term ensures the convergence of the series. Furthermore, the derivative of Weierstrass elliptic function is an odd function satisfying the differential equation  
\[
(\wp^{\prime}(z))^{2}=4(\wp(z))^{3}-g_{2}\wp(z)-g_{3},
\]
where 
$g_{2} := \sum_{\ell \neq 0}60\frac{1}{\left(\ell\right)^{4}}$
and
$g_{3} := \sum_{\ell \neq 0}140\frac{1}{\left(\ell \right)^{6}}$. 
This differential equation can be used to compute higher-order derivatives of $\wp$. 
We obtain 
\begin{align*}
\wp^{(2)}(z) &= 6 \wp^2(z) - \frac{g_2}{2}  
\end{align*}
and 
$$
\wp^{(n+2)}(z) = 6 \sum_{k=0}^n  \binom{n}{k} \wp^{(n-k)}(z) \wp^{(k)}(z), 
\qquad \qquad n \geq 1. 
$$

The Weierstrass zeta function is defined by 
\begin{equation} \label{e:zeta}
\zeta(z)=\frac{1}{z}+\sum_{\ell \neq 0}\left\{ \frac{1}{z-\ell}+\frac{1}{\ell}+\frac{z}{\ell^{2}} \right\}=\frac{1}{z}+\sum_{\ell \neq 0} \frac{z^3}{\ell^2 (z^2-\ell^2)}   
\end{equation}
and satisfies
\[
\frac{d\zeta}{dz}=- \wp (z).
\]
It has a Laurent expansion near $z=0$
\[
\zeta(z)=\frac{1}{z}-\sum_{k=2}^{\infty} \gamma_{2k} z^{2k-1},   
\]
where 
$\gamma_{2k} = \sum_{\ell \neq 0} \frac{1}{\ell^{2k}}$,  $k\ge2$.
In contrast to $\wp (z)$, the function $\zeta(z)$ does not possess the double-periodic property. Instead, it satisfies the quasi-periodic condition:
\[
\zeta(z+2\omega_{1})=\zeta(z)+2\eta_{1}
\]
\[
\zeta(z+2\omega_{2})=\zeta(z)+2\eta_{2}
\]
where $\eta_{1}=\zeta(\omega_{1})$ and $\eta_{2}=\zeta(\omega_{2}).$
The values $\eta_{1}$, $\eta_{2}$, $\omega_{1}$, $\omega_{2}$
are not independent but related by the Legendre identity 
\[
\eta_{1}\omega_{2}-\eta_{2}\omega_{1}= \frac{\pi \imath}{2}.
\]
The $\zeta$ function can be modified so that it is periodic,  
$$
\hat{\zeta}(z) = \zeta(z)- \gamma_2 z- \frac{\pi}{A} z^{*}
$$ 
where 
$A$ is the area of the fundamental cell of the lattice and 
$\gamma_2 $ is given by a Eisenstein summation and satisfies 
$\zeta(\omega_i) \equiv \eta_i = \gamma_2 \omega_i +\frac{\pi \omega_i^{*}}{A}$, $i=1,2$ \cite{Haldane_2018}.  
Note that since $\hat{\zeta}$ depends on $z^*$, it is no longer meromorphic. 

Finally, the Weierstrass sigma function is defined by 
$$
\sigma(z,\omega) = \lim_{\varepsilon \to 0} \varepsilon \exp \left( \int_\varepsilon^z \zeta(w,\omega) dw \right), 
$$
which is an odd, non-doubly-periodic, holomorphic function with simple zeros at the lattice points. 
 It satisfies 
\begin{equation} \label{e:zeta-omega-relationship}
\zeta(z,\omega) = \frac{\sigma'(z,\omega)}{ \sigma(z,\omega)}. 
\end{equation}
As for the zeta function, the sigma function can be modified as in \eqref{e:sigma_hat}. When defined this way, its modulus has the lattice periodicity \cite{Haldane_2018}.

\section{Proof of Theorems \ref{t:LogConjThmTorus} and \ref{t:HarmonicPerDecomp}}
\label{s:LogConjThmProof}
\begin{proof}[Proof of \cref{t:LogConjThmTorus}.]
The first part of the proof closely follows the proof of  S. Axler for the Logarithmic Conjugation Theorem \cite{Axler_1986}. 
Define $h\colon \Omega \to \mathbb C$ by 
$$
h(z) := u_x(z) - \imath u_y(z). 
$$
The Cauchy-Riemann equations can be used to check that $h$ is analytic on $\Omega$. 
For each $j\in [b]$, let $\Gamma_j$ be a closed curve in $\Omega$ that circles $K_j$ once and no other $K_k$, $k\neq j$. Define
$$
c_j := \frac{1}{2 \pi \imath } \oint_{\Gamma_j} h(w) \ dw. 
$$
We see that $\Im c_j 
= - \frac{1}{2\pi} \Re \oint_{\Gamma_j} h(w) \ dw 
= - \frac{1}{2\pi} \Re \oint_{\Gamma_j} u_x(w) dx + u_y(w) dy = 0$, so $c_j$ is a real number for each $j\in [b]$. 
Since $u$ is doubly-periodic, so is $h$, and by the Cauchy Integral Theorem \cite[Thm.1]{Datar}, we have that 
\begin{equation} \label{e:sum-c=0}
\sum_{j \in [b]} c_j = 0.
\end{equation}

We consider $h$ to be a function on $\mathcal P$, which we still denote by $h$. 
Fix a point $z_0 \in \mathcal P$, and define $f\colon \mathcal P \to \mathbb C$ by
$$
f(z) := \int_{z_0}^z h(w) - \sum_{j\in [b]} c_j  \zeta( w - a_j, \omega) \ dw, 
$$
where the integral is taken over any path in $\mathcal P$ from $z_0$ to $z$
and $\zeta$ is the Weierstrass zeta function as in \eqref{e:zeta}. 
To show that $f$ is well-defined, we check that the above integral is independent of the path from $z_0$ to $z$. Take two paths from $z_0$ to $z$ and reverse the direction of transversal in one to form a closed curve. Thus, we need only show that 
$$
\frac{1}{2\pi \imath } \oint_{\gamma} h(z) dw =  
\frac{1}{2\pi \imath } \sum_{j\in [b]} c_j \oint_\gamma  \zeta( w - a_j, \omega) \ dw
$$
for any closed curve $\gamma$ $\Omega$. By the Cauchy Integral Theorem and the definition of $c_j$, the left hand side is given by $\sum_{j \in [b]} c_j I_j(\gamma) $, where $I_j(\gamma)$ denotes the winding number of $\gamma$ about $K_j$. 
Using that the Laurent expansion for 
$\zeta(z,\omega)$, which has a single pole of order one, by the Cauchy Integral Theorem,  the right hand side is also seen to be equal to $\sum_{j \in [b]} c_j I_j(\gamma)$, as desired.  The function $f(z)$ is analytic on $\mathcal P$ and we compute the derivative 
\begin{equation} \label{e:f'}
f'(z) = h(z) - \sum_{j\in [b]} c_j \zeta( z - a_j, \omega).  
\end{equation}

Now define 
\begin{equation} \label{e:q}
q(z) := \Re f(z) + \sum_{j \in [b]}  c_j \log | \sigma \left(z - a_j, \omega \right) |. 
\end{equation}
We claim that $u_x(z) = q_x(z)$ and $u_y(z) = q_y(z)$, so that, after adding a constant to $f$, we obtain $u(z) = q(z)$, $z \in \mathcal P$. Using \eqref{e:zeta-omega-relationship},  we compute
\begin{align*}
q_x(z) &= \Re f'(z) + \sum_{j \in [b]}  c_j \Re \zeta(z - a_j, \omega) 
= \Re h(z) 
= u_x. 
\end{align*}
and
\begin{align*}
q_y(z) &= \Re  \left( \imath f'(z) \right)  + \sum_{j \in [b]}  c_j \Re \left( \imath \zeta(z - a_j, \omega) \right) 
= \Re \left( \imath h(z) \right)  
= u_y. 
\end{align*}

We have established that $u(z) = q(z)$ up to a constant, $z \in \mathcal P$ and it remains to show that we can rewrite $q(z)$ in \eqref{e:q} so that the two terms on the right hand side are each doubly-periodic, so can be thought of as functions on $\Omega$. In \eqref{e:q}, the second term on the right hand side is not doubly-periodic since $\sigma$ is not doubly-periodic. By \eqref{e:sigma_hat}, this term can be rewritten 
$$
\sum_{j \in [b]}  c_j   \log | \sigma \left(z - a_j\right) | = 
\sum_{j \in [b]}  c_j  \log | \hat \sigma \left(z - a_j\right) | 
+ \Re g(z)
$$
where 
\begin{align*}
g(z) 
&= \frac{1}{2} \sum_{j \in [b]}  c_j \left( 
 \gamma_2 (z-a_j)^2  + \pi |z-a_j|^2/A \right) \\
&= \alpha z + \beta z^* + \gamma, 
\end{align*}
where $\alpha$, $\beta$, and $\gamma$ are constants and we have used \eqref{e:sum-c=0} to drop the quadratic terms. 

From \eqref{e:f'}, $f'$ is doubly-periodic since $h$ is doubly-periodic and $\sum_{j \in [b]} c_j = 0$. 
There exists $\alpha_1, \alpha_2$ such that for all admissible $z$
$$
\left\{
    \begin{array}{l}
f(z+2\omega_{1})  = f(z) + \alpha_1 \\
f(z+2\omega_{2})  = f(z) + \alpha_2. 
\end{array}
\right.
$$
Let us introduce $(\mu_1, \mu_2)$ the unique solution of 
$$
\left\{
    \begin{array}{l}
\omega_1 \mu_1 + \omega_1^* \mu_ 2 = - \alpha_1 / 2.\\
\omega_2 \mu_1 + \omega_2^* \mu_ 2 = - \alpha_2 /2.
\end{array}
\right.
$$
Notice that previous system is non-singular since the determinant is proportional to the area of the fundamental domain, which is nonzero. Moreover, a straightforward computation shows that 
$$f(z) + \mu_1 z + \mu_2 z^*$$
is a doubly-periodic function. Thus,  for a suitable $\mu$, 
$\hat{f}(z) = f(z) + \mu z$, is an analytic function and $\Re \hat{f}(z) $ is also doubly-periodic. Note that $\Im \hat f(z)$ is not necessarily doubly-periodic and $\hat f '$ is elliptic. 

Summarizing our results, we have established that there exists $\hat{f}$ analytic with doubly-periodic real part  and $(\nu, \xi) \in \mathbb{C}^2$ such that 
$$u = \Re \hat{f} +\sum_{j \in [b]}  c_j  \log | \hat \sigma \left(z - a_j\right) | + \Re{(\nu z + \xi z^*)}$$
Observing that both the left hand side and the two first terms of the right hand side are doubly-periodic, we obtain $\nu = \xi = 0$, which concludes the proof.
\end{proof}

\begin{proof}[Proof of \cref{t:HarmonicPerDecomp}.]
Let $\hat f'$ be the elliptic function  from \cref{t:LogConjThmTorus} associated with the harmonic function $u$. Using a representation of elliptic functions (see, e.g., \cite[p.450]{whittaker1920course} or \cite[p.23]{pastras2020weierstrass}), we may write 
$$
\hat f'(z) = \tau + \sum_{j\in [b]} \left( \alpha_{j}  \zeta(z - a_j) + \sum_{k\ge0} \beta_{j,k}  \wp^{(k)}(z-a_j)  \right), 
$$
where $\tau \in \mathbb C$, $\alpha_j \in \mathbb C$, and $\beta_{j,k} \in \mathbb C$ are constants. 
Consequently, there exists a constant $\rho \in \mathbb C$ such that 
$$
\hat f(z) = \tau z + \rho + \sum_{j\in [b]} \left( \alpha_{j}  \log \sigma(z - a_j) + \beta_{j,0} \zeta(z - a_j) +  \sum_{k\ge1} \beta_{j,k}  \wp^{(k- 1)}(z-a_j)  \right). 
$$
Introducing the periodic modifications $\hat{\zeta}$ and $\log |\hat{\sigma}|$ of $\zeta$ and $\log |\sigma|$ functions respectively, we obtain that there exists real coefficients $C,a_{j,k}, b_{j,k}$ such that
{\small \begin{equation*} 
\begin{split}
u(z) =  C +  \sum_{j\in [b]} & \Big[ a_{j,-1} \Re \hat{\zeta}(z - a_j) + 
                                    b_{j,-1} \Im \hat{\zeta}(z - a_j) +
            \sum_{k\ge0} a_{j,k}  \Re \wp^{(k)}(z-a_j) + 
             \sum_{k\ge0} b_{j,k}  \Im \wp^{(k)}(z-a_j) \\
           & + c_j \log | \hat \sigma \left( z - a_j \right) | \Big] + g(x,y)
\end{split}
\end{equation*}}for some affine function $g$. By periodicity of all other terms, the function $g$ has also to be doubly-periodic, so must be identically equal to zero. 
Finally, $\sum_{j\in [b]} c_j = 0$ is deduced from the harmonicity of all the terms, except the $\log|\hat \sigma|$ terms which have a constant Laplacian.
\end{proof}

\section{Computational method and experiments}
\label{s:examples}
Here we develop a computational method based on a series solution of the form \eqref{e:HarmonicRep} to solve the 
Laplace problem \eqref{e:Laplace} and the 
Steklov eigenvalue problem \eqref{e:Steklov}.

\subsection{Computational Method}
Let $\Omega$ be a finitely-connected torus as in \eqref{e:Omega}. For simplicity, we will take each region $K_j$, $j\in [b]$ to be a closed disk, $K_j = \overline{B}(a_j, r_j)$, that is centered at the point $a_j$ and has radius $r_j$. The centers and radii are chosen such that $K_i \cap K_j = \varnothing$ for $i\neq j$. 
Based on \cref{t:HarmonicPerDecomp}, we consider a series solution of the form \eqref{e:HarmonicRep}, where the sums on $k$ are truncated at $k=k_{\max}$. 
We collect the (real) coefficients in the series solution into a vector $v = [C, (a_{j,k}), (b_{j,k}), (c_j) ] \in \mathbb R^m$, where $m=1+2 b (k_{\max}+2) + (b-1)$. For each coefficient, $v_i$, we let $\phi_i$, $i \in [m]$ denote the corresponding basis function (e.g., the real part of a Weierstrass $\wp$ function), so that 
\begin{equation}
\label{e:truncatedSeriesSolution}
u(z) = \sum_{i \in [m]} v_i \phi_i(z).
\end{equation}
On each boundary component $\partial K_j$, we uniformly sample points with respect to arclength and denote the collection of all sampled points in the union of the boundary components by $(p_\ell)_{ \ell \in [S]}$. 
In the experiments below, we report the value of $m$ and take $S = 3 m$.  
Define the matrices $A,B \in \mathbb R^{S\times m}$ by  
\begin{align*}
A_{\ell,i} &= \frac{\partial \phi_i}{\partial n}(p_l)  \\
 B_{\ell,i} &= \phi_i(p_l). 
\end{align*}
Details about the computation of the normal derivatives of basis functions are given in \cref{s:NormalDerivatives}.

\begin{figure}
  \centering
  \includegraphics[width=0.45\textwidth]{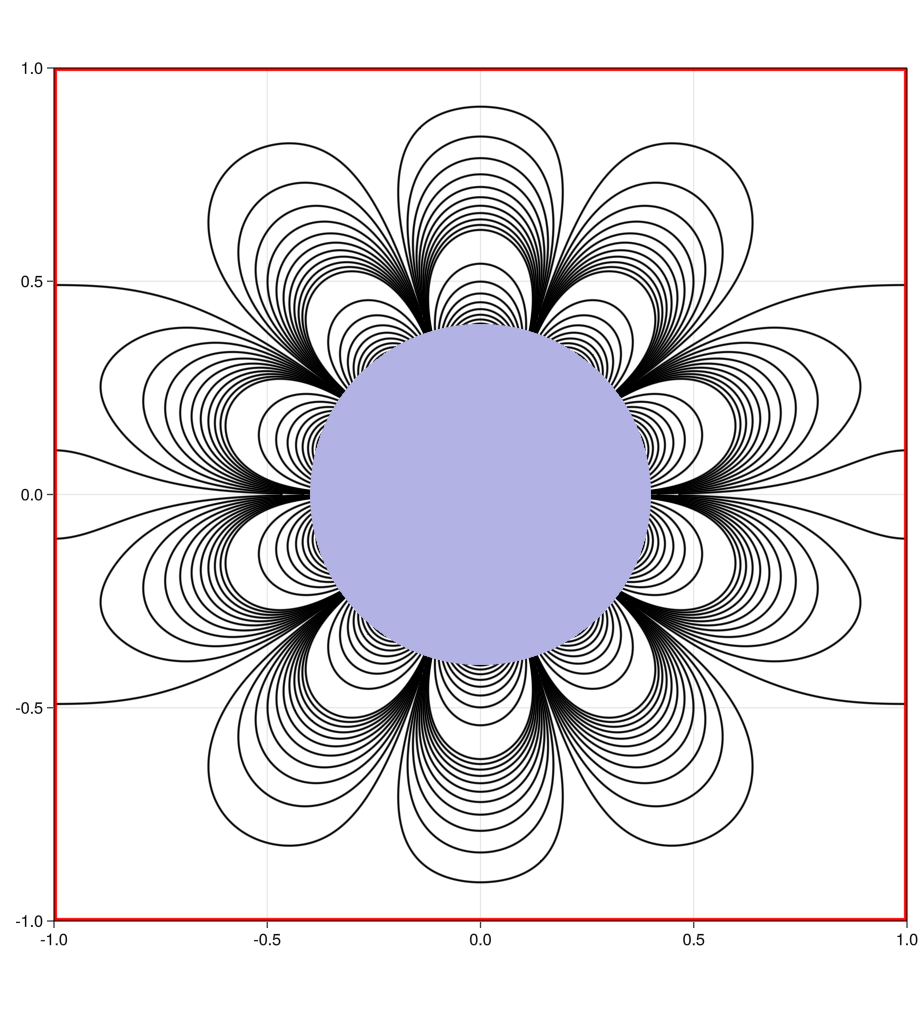}
  \includegraphics[width=0.45\textwidth]{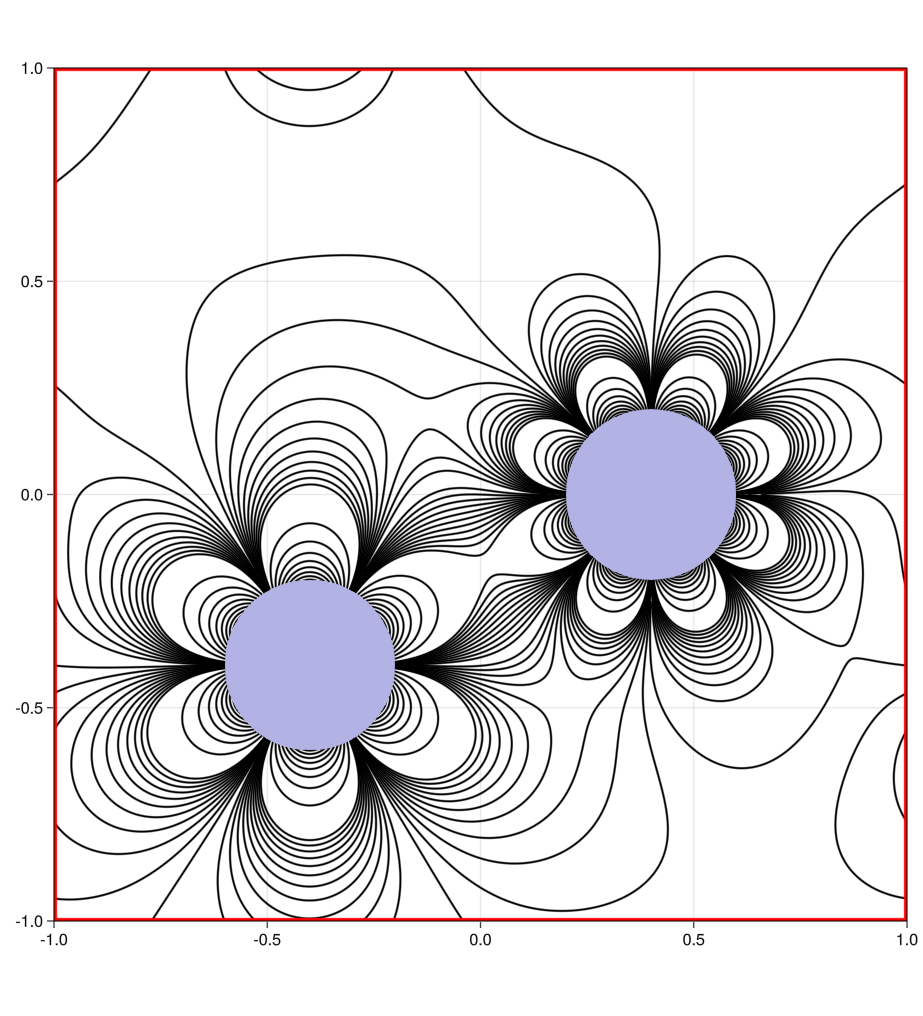} \\
  \includegraphics[width=0.45\textwidth]{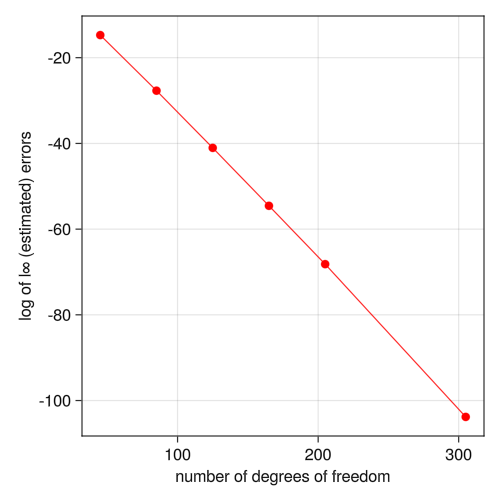}
  \includegraphics[width=0.45\textwidth]{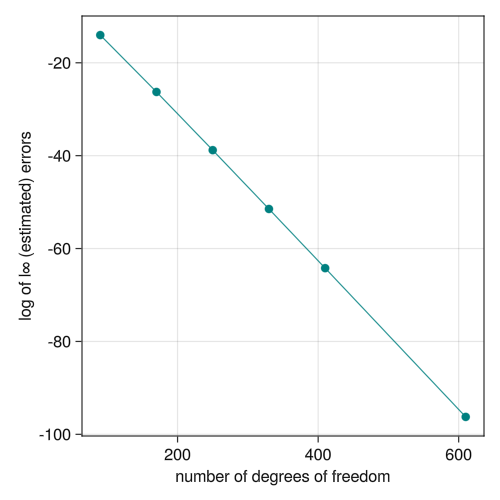} 
  \caption{ {\bf (Upper panels)} Approximate solution to the Laplace problem in a square torus with one and two circular holes. {\bf (Lower panels)} Spectral convergence is observed for each of the two geometries. See \cref{s:Laplace}. }
\label{f:LaplaceSquare}
\end{figure}

\begin{figure}
  \centering
  \includegraphics[width=0.45\textwidth,trim={0 8cm 0 8cm},clip]{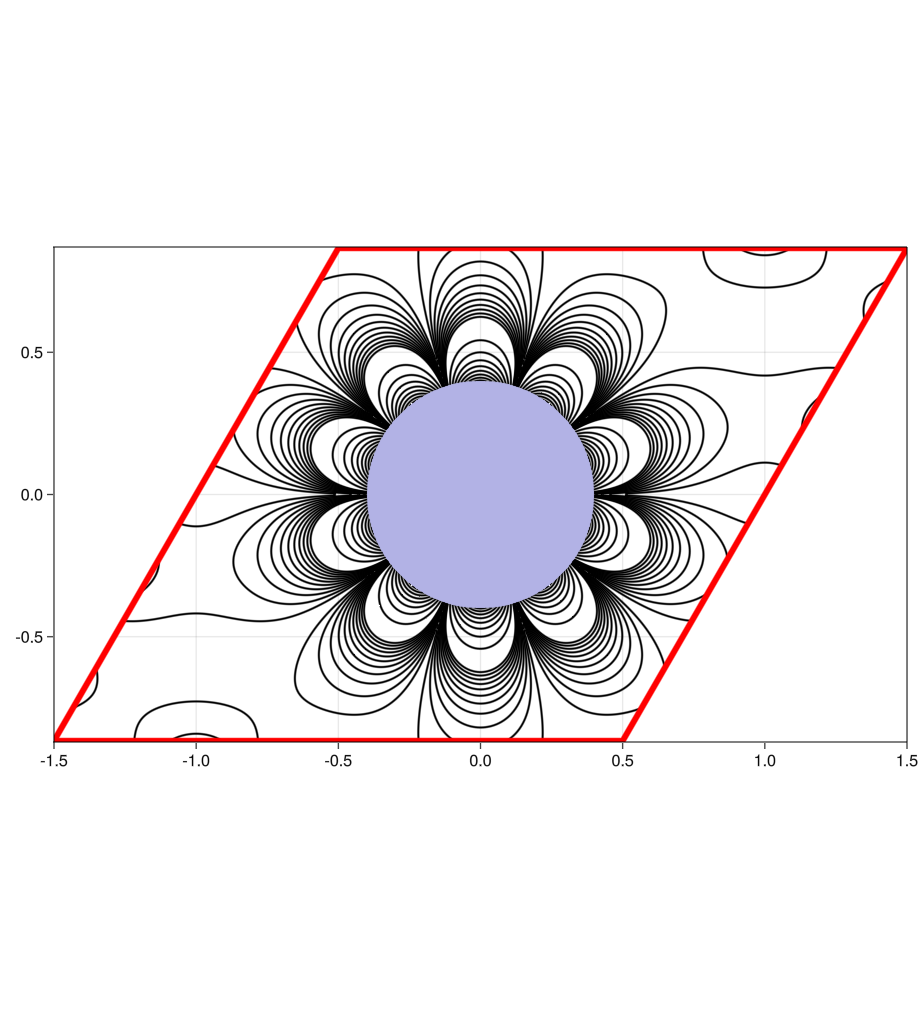}
  \includegraphics[width=0.45\textwidth,trim={0 8cm 0 8cm},clip]{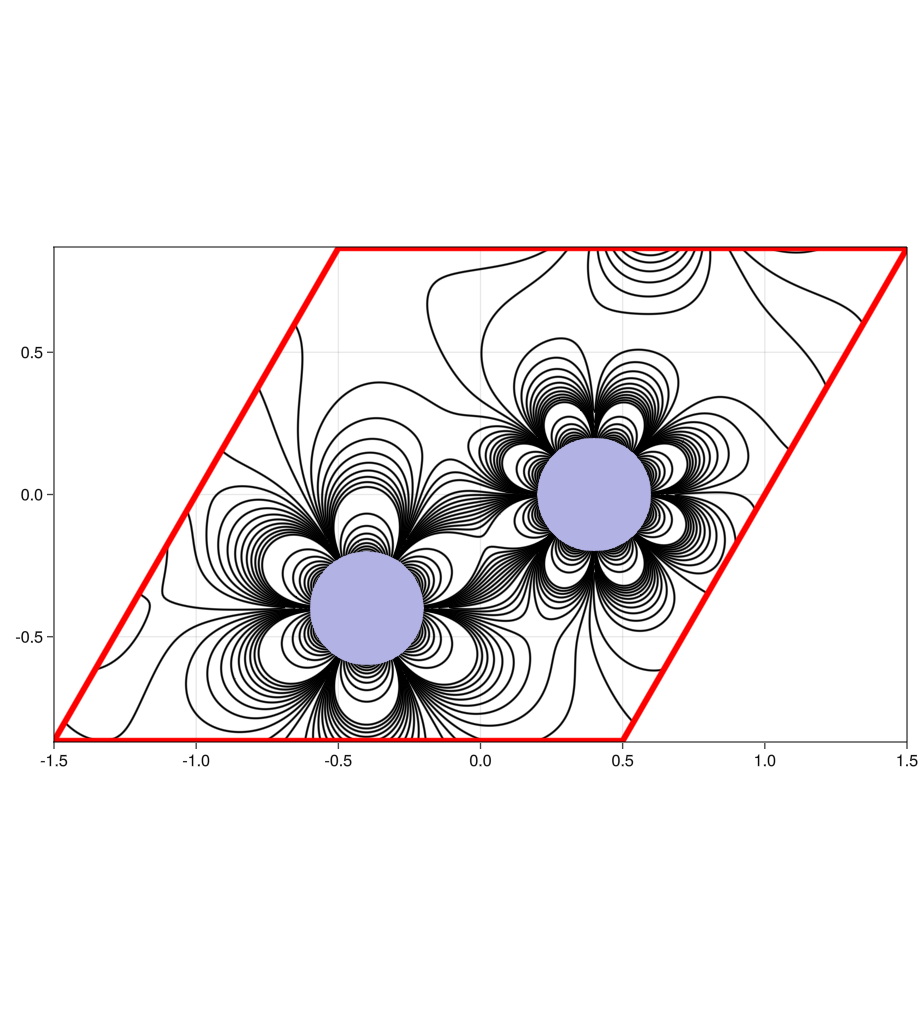}
  \caption{Approximate solution to the Laplace problem in an equilateral torus with one and two circular holes. See \cref{s:Laplace}. }
  \label{f:LaplaceEquilateralTorus}
\end{figure}

\subsection{Laplace problem} \label{s:Laplace}
We solve the Laplace problem \eqref{e:Laplace}, with boundary data $f(x)$, $x \in \partial \Omega$ as follows. Define the vector $b \in \mathbb R^S$ by $b_\ell = f(p_\ell)$. The least-squares solution is found by solving the normal equations \begin{equation}
\label{e:normalEq}
B^t B v = B^t b.
\end{equation}
The solution $v$ then is used with the expansion in \eqref{e:truncatedSeriesSolution} as an approximate solution of \eqref{e:Laplace}. By the maximum principle, the accuracy of the solution can be computed by looking at the error on the boundary, 
$\sup_{x\in \partial \Omega} |u(x) - f(x)|$.

We implement the numerical method in Julia using arbitrary precision provided by the packages \emph{GenericLinearAlgebra.jl} and \emph{ArbNumerics.jl} (a wrapper of the \emph{Arb} C library). All computational experiments were performed with a precision of $2^{10}$ bits which corresponds to a machine epsilon approximately equal to $10^{-300}$. 

We first consider a finitely-connected square torus with 
half-periods $(\omega_1,\omega_2) = (1,i)$. The complement is taken to be $b=1$ disks with
$K_1 = B(a_1,r)$ with $a_1 = 0$ and $r=0.4$. 
We take $f(\theta) = \sin(5\theta)$, where $\theta$ is the polar angle centered at $a_1$. 
The resulting solution is plotted in the top left panel of \cref{f:LaplaceSquare}. 
Using the maximum principle to bound the $L^\infty(\Omega)$ error of the solution, we estimate $\| u(x) - f(x) \|_{L^\infty(\partial \Omega)}$ in the bottom left panel of \cref{f:LaplaceSquare} for increasing number of 
degrees of freedom,  $m$. This estimate is based on the maximum value obtained at the sampled points, after doubling the number of sampled points. 
Spectral convergence is observed. 
With $k_{\max} = 150$ ($m = 305$ degrees of freedom), the solution has error less $10^{-100}$ corresponding to at least 100 digits of accuracy.  

Next, we again consider a finitely-connected square torus with 
half-periods $(\omega_1,\omega_2) = (1,i)$. The complement is taken to be  
$b=2$ disks with 
$K_i = B(a_i,r_i)$, $i=1,2$ with $a_1 = 0.4$, $a_2 = -0.4-0.4i$ and $r_1 = r_2 = 0.2$. 
On each boundary, we set 
$f(\theta) = \sin(4\theta)$ for circle $i=1$ and 
$f(\theta) = \sin(3\theta)$ for circle $i=2$. 
The resulting solution is plotted in the top right panel of \cref{f:LaplaceSquare}. 
In the bottom right panel of \cref{f:LaplaceSquare}, we can see that the solutions have similar error as the previous example, albeit using more degrees of freedom. 

Next, we consider a finitely-connected equilateral torus with 
half-periods $(\omega_1,\omega_2) = \left(1,\frac{1}{2} + \frac{\sqrt{3}}{2} \imath \right)$. The complement is taken to be the same sets as above with one and two circular holes. We plot the results in \cref{f:LaplaceEquilateralTorus}. The solutions have similar error to the previous examples.

\begin{figure}
  \centering
  \includegraphics[width=0.45\textwidth,trim={0 8cm 0 8cm},clip]{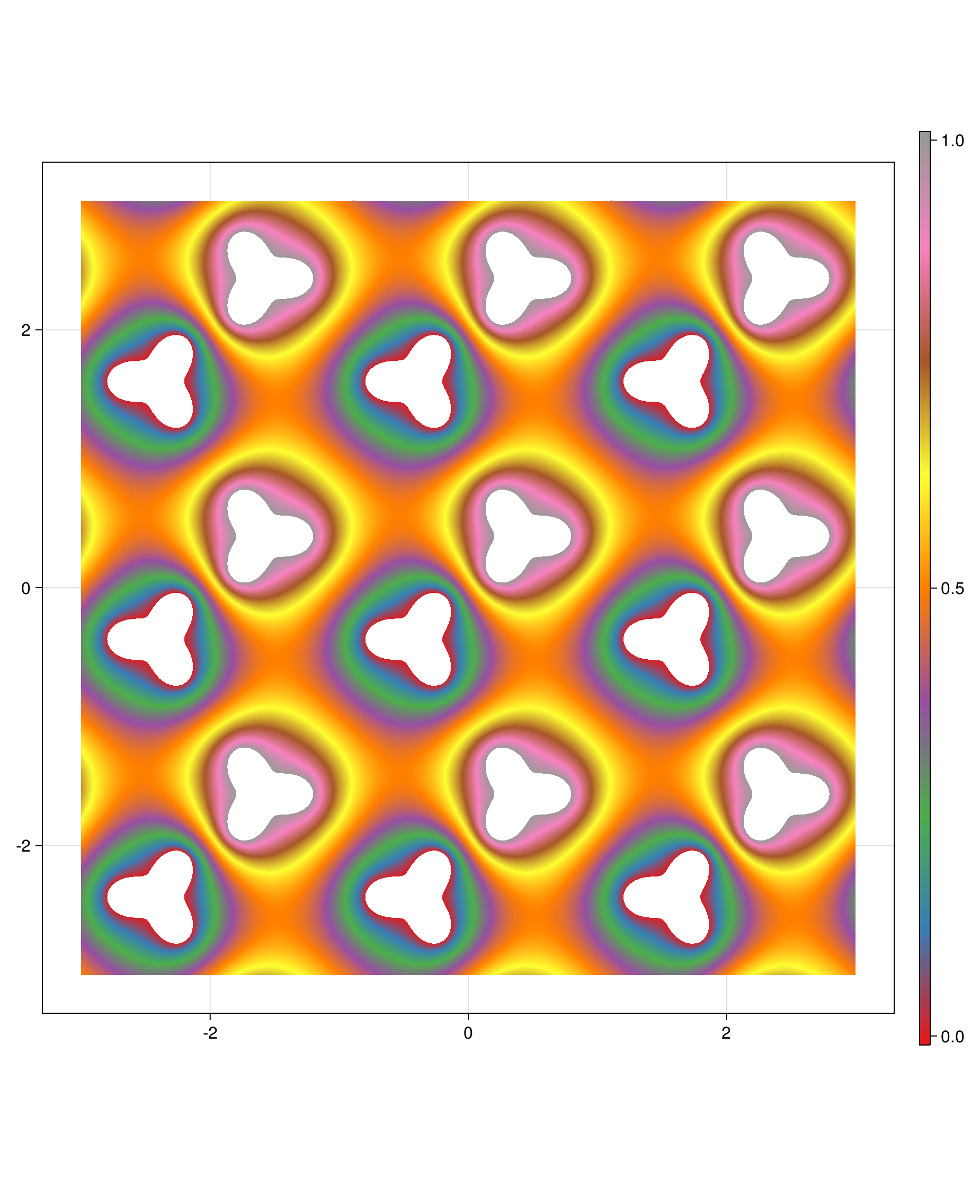}
  \includegraphics[width=0.45\textwidth,clip]{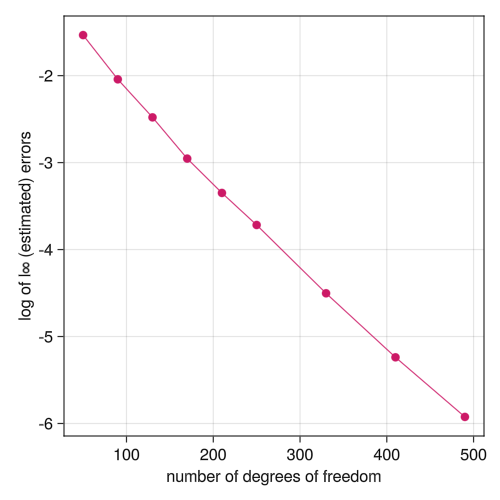}
  \caption{Approximate solution to the Laplace problem in a square torus with two-non convex holes.}
  \label{f:LaplaceFlowers}
\end{figure}

In \cref{f:LaplaceFlowers}, we provide an approximate solution to the Laplace problem for two non-convex holes in a square torus. The polar parametrizations of the boundaries of the two holes $K_1$ and $K_2$ are given by $A_i + \rho(\theta + \theta_i) (\cos(\theta ),\sin(\theta))$ where
$$
\rho(\theta) = \frac{3}{10}+\frac{1}{10}\cos(3\theta),
$$
$A_1 = (0.4, 0.4)$, $A_2 = -A_1$, $\theta_1 = 0$, and $\theta_2 = \frac{\pi}{3}$. We impose the Dirichlet condition $0$ on the first boundary component and $1$ on the second. The sampled points are obtained using the previous parametrization together with a uniform sampling of the angles. As previously, in the right panel of \cref{f:LaplaceFlowers}, we observe exponential convergence but notice that the obtained accuracy is significantly lower than in previous cases with the same number of degrees of freedom. 

Finally, we consider the Laplace equation on a finitely-connected square torus with 25 disks removed. Dirichlet boundary conditions equal to 0 or 1 are imposed on the boundary of each disk. The results are plotted in \cref{f:multiholes}. The solution has error less than $10^{-16}$. 

\begin{figure}[t]
  \centering
  \includegraphics[width=0.45\textwidth,clip]{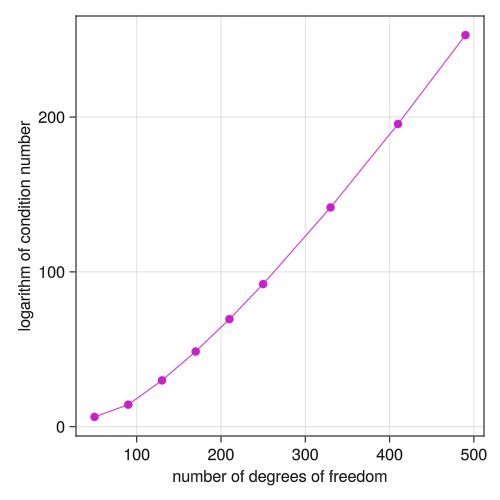}
  \includegraphics[width=0.45\textwidth,clip]{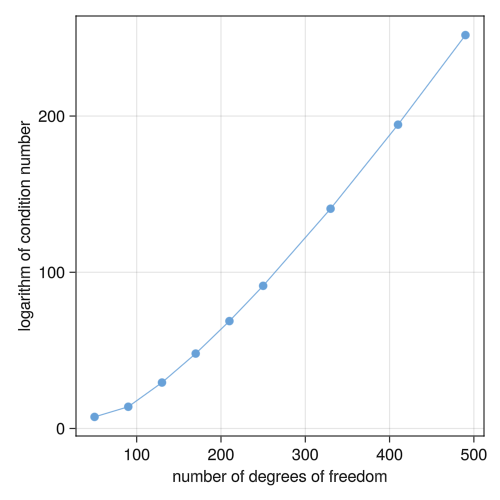}
  \caption{Exponential growth of condition numbers of matrices $B^t B$ (left) and $B^t A$ (right) with respect to the number of degrees of freedom.}
  \label{f:condnumbers}
\end{figure}

\subsection{Steklov eigenvalue problem} \label{s:Steklov}
To solve the Steklov eigenvalue problem \eqref{e:Steklov}, we consider a generalized eigenvalue problem
\begin{equation} \label{e:mpseig}
A  \ v =  \sigma  B \  v. 
\end{equation} 
We can approximate solutions to this eigenvalue problem by multiplying both sides by $B^t$ and considering the non-symmetric generalized eigenvalue problem, $B^t A v = \sigma B^t B v$. 
For $k_{max} \leq 50$, this formulation leads to  exponentially converging eigenvalue approximations, as expected. As it has been observed by several authors \cite{betcke2005reviving, betcke2008generalized,Fox_1967}, this formulation with a larger number of degrees of freedom may produce ill-conditioned matrices. 
To illustrate this, in \cref{f:condnumbers}, for the example considered above with two non-convex holes (see \cref{f:LaplaceFlowers}), we plot the condition number of $B^tB$ and $B^tA$ as the number of degrees of freedom varies. 
To overcome this difficulty and avoid spurious modes, we followed the SVD approach described in  \cite{betcke2005reviving}: for a (small) set of randomly sampled interior points $(q_r)_{r \in [R]}$ we consider the evaluation matrix  $C \in \mathbb R^{R\times m}$ 
 \begin{equation}
  C_{r,i} = \phi_i(q_r) 
\end{equation} 
In all our experiments we set $R = 50$. We define $s(\sigma)$ to be the smallest (always non-negative) eigenvalue of the generalized eigenvalue problem 
\begin{equation}
  D(\sigma) x(\sigma) = s(\sigma) C^t C x(\sigma)
\end{equation}
where $D(\sigma) = (A - \sigma B)^t (A - \sigma B)$. From a computational point of view, $s(\sigma)$ can be efficiently evaluated using a standard power method or an orthogonal subspace approach if the multiplicity is suspected to be greater than one. Local minimizers of 
 $s(\sigma)$ provide stable approximations of Steklov eigenvalues. To identify numerically these local extrema, we used the golden section algorithm. 

To bound the error in the eigenvalues, we use the following a posteriori estimate for Steklov eigenvalues in  \cite{Bogosel_2016}, 
which extends previous estimates for Laplace-Dirichlet eigenvalues \cite{Fox_1967,Moler_1968}. 
\begin{prop}[\cite{Bogosel_2016}]
  \label{p:steklovbound}
Consider $\Omega$ a bounded open regular domain, and suppose that $u_\varepsilon$ solve the following approximate eigenvalue problem 
\begin{align*}
& - \Delta u_\varepsilon = 0 && \textrm{in } \Omega \\
& \partial_n u_\varepsilon = \sigma_\varepsilon u_\varepsilon + f_\varepsilon  && \textrm{on } \partial \Omega. 
\end{align*}
Then if $\| f_\varepsilon \|_{L^2(\partial \Omega)}$ is small, there exists a constant $C$, depending only on $\Omega$, and a Steklov eigenvalue $\sigma_k$ satisfying 
$$
\frac{|\sigma_\varepsilon - \sigma_k|}{ \sigma_k} \leq  C \| f_\varepsilon \|_{L^2(\partial \Omega)}. 
$$
\end{prop}

\begin{figure}
  \centering
  \includegraphics[width=0.28\textwidth]{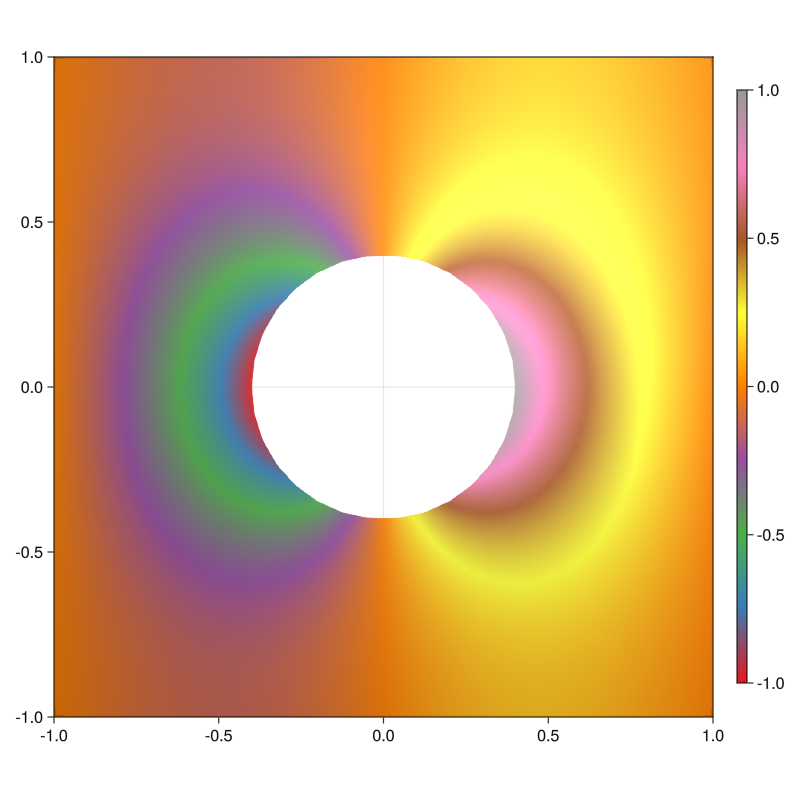}
  \includegraphics[width=0.28\textwidth]{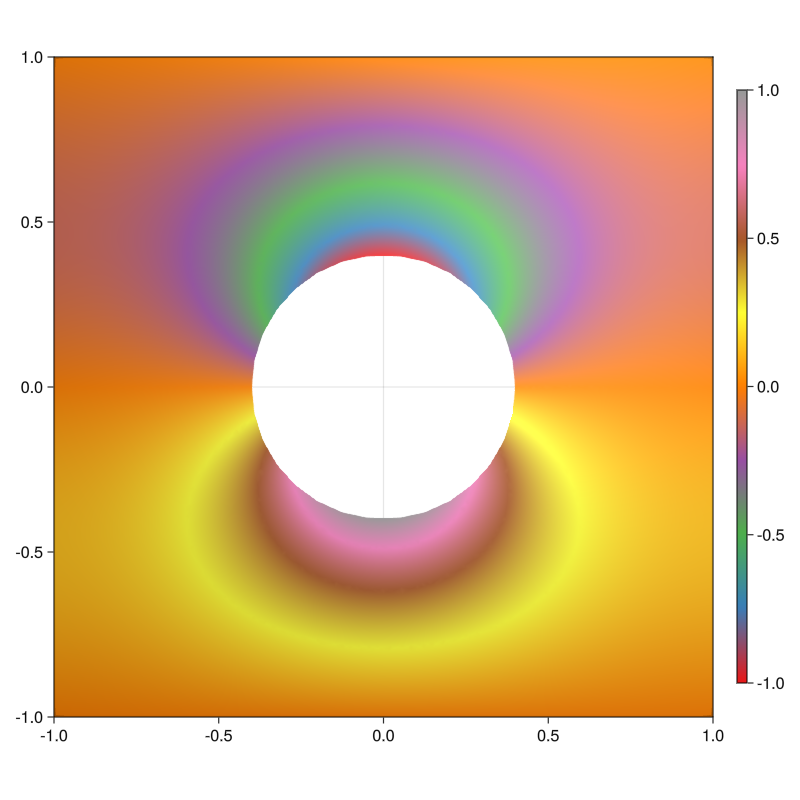}
  \includegraphics[width=0.28\textwidth]{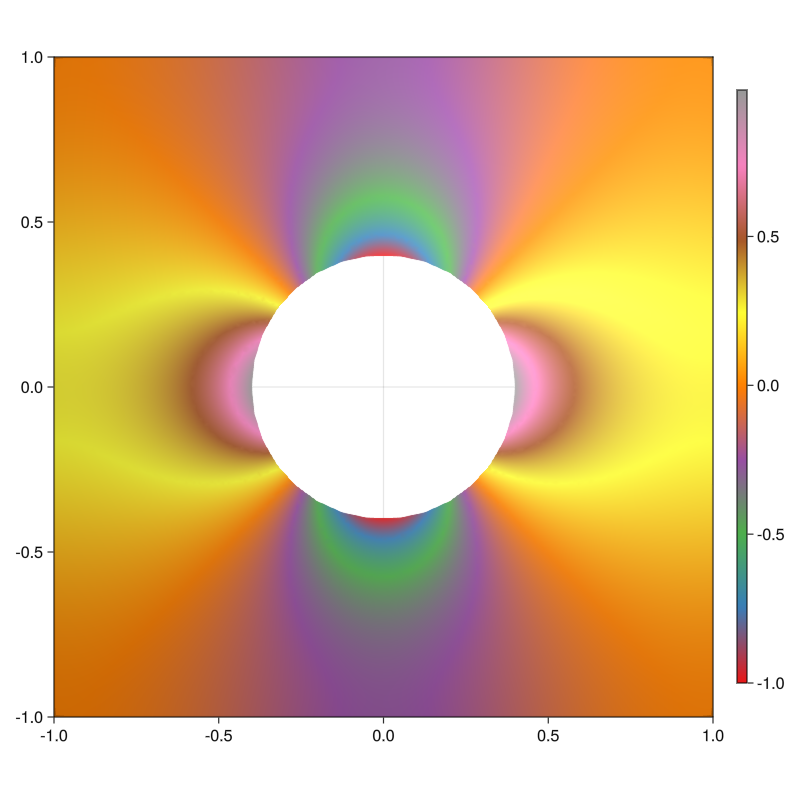}
  \includegraphics[width=0.28\textwidth]{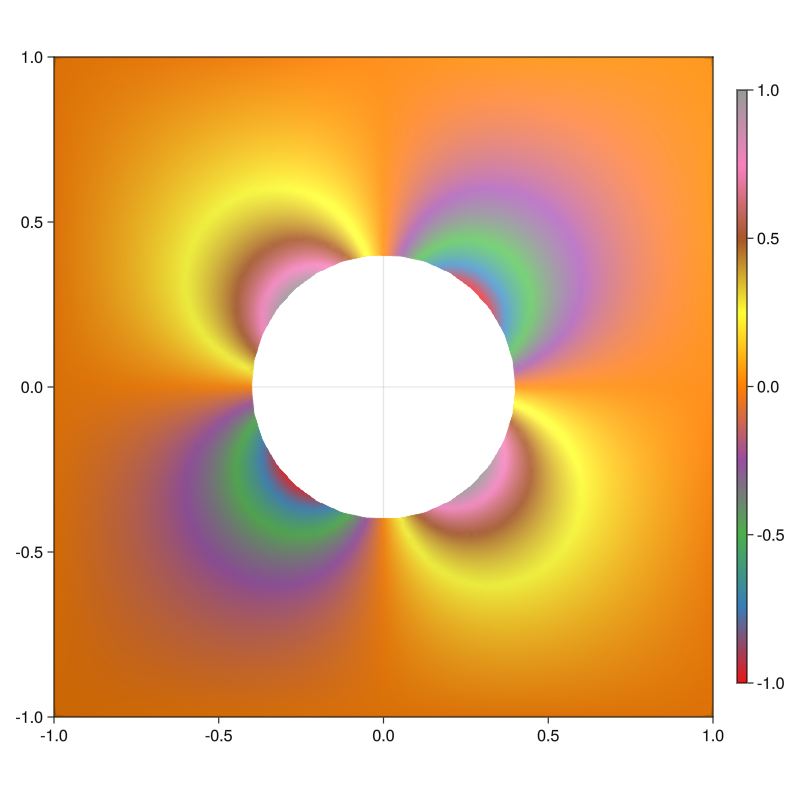}
  \includegraphics[width=0.28\textwidth]{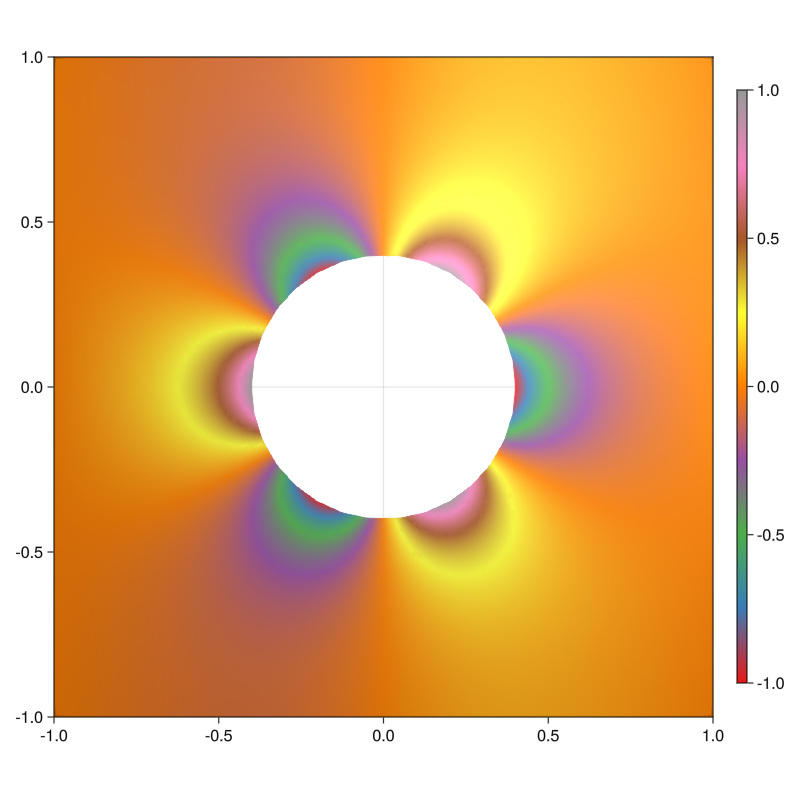}
  \includegraphics[width=0.28\textwidth]{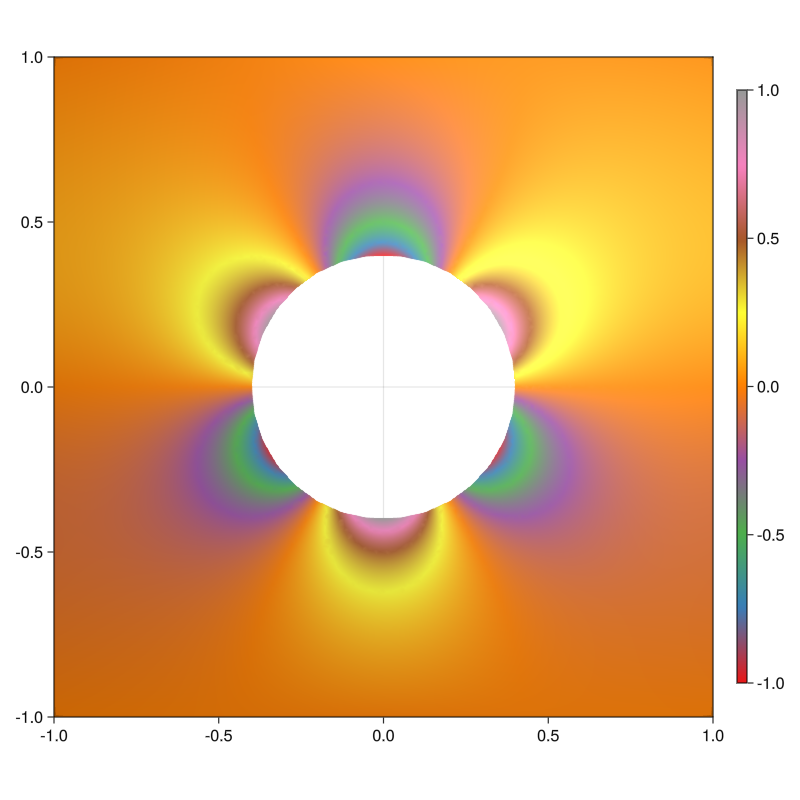}
  \caption{Approximate Steklov eigenfunctions of indices $2$ to $7$ on a punctured square torus with one hole. See \cref{s:Steklov}.}
  \label{f:SteklovEigSquare1} 
\end{figure}

\begin{figure}[p]
  \centering
  \includegraphics[width=0.28\textwidth]{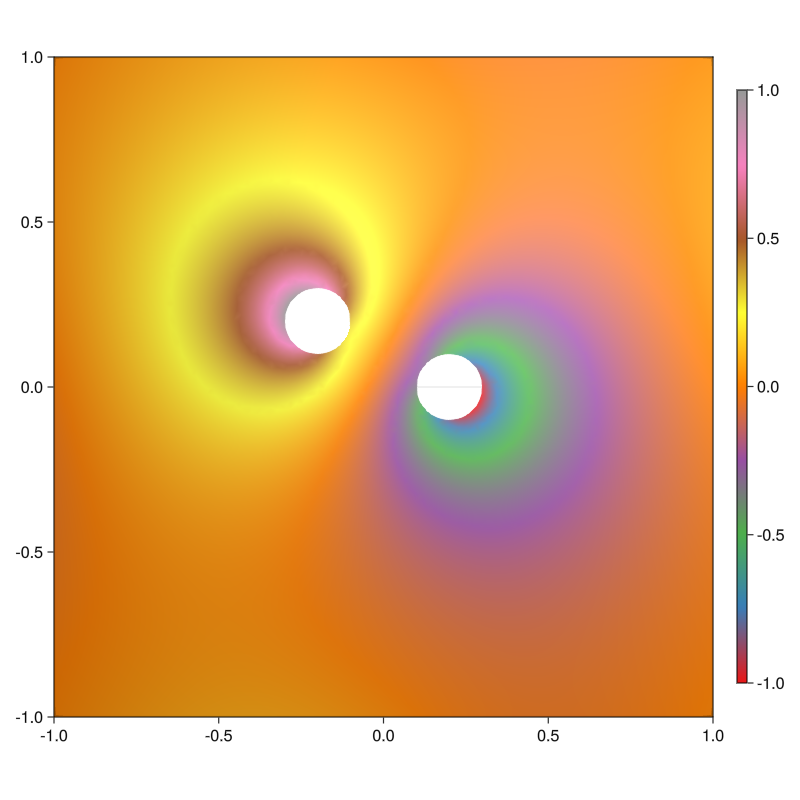}
  \includegraphics[width=0.28\textwidth]{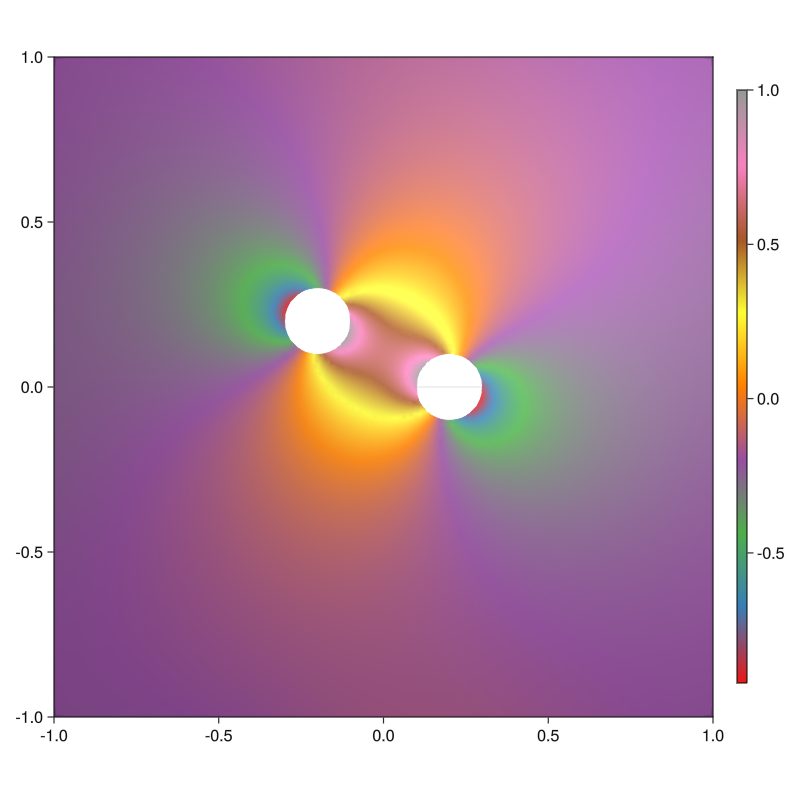}
  \includegraphics[width=0.28\textwidth]{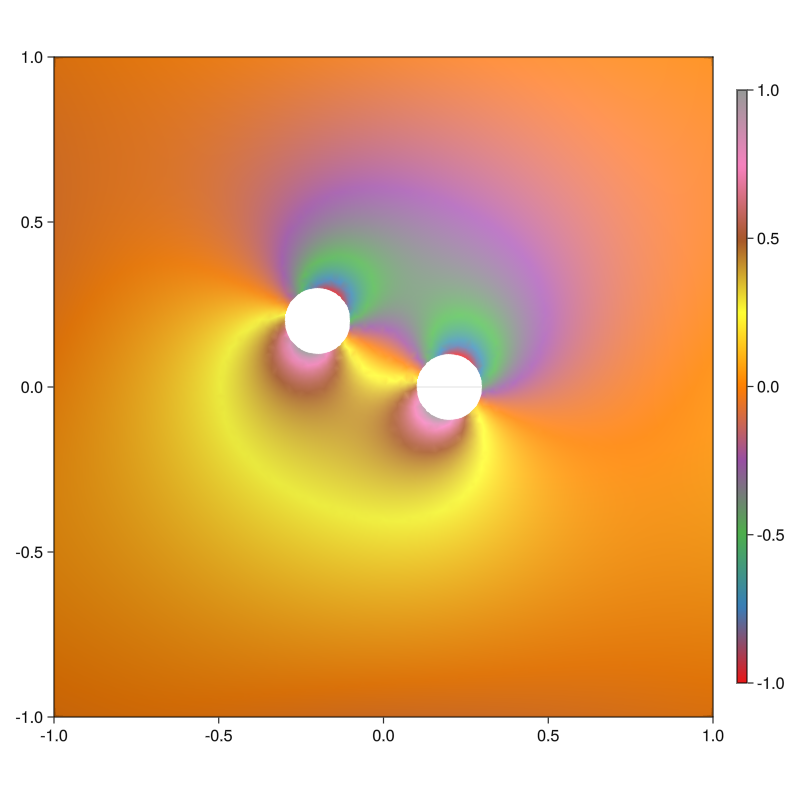}
  \includegraphics[width=0.28\textwidth]{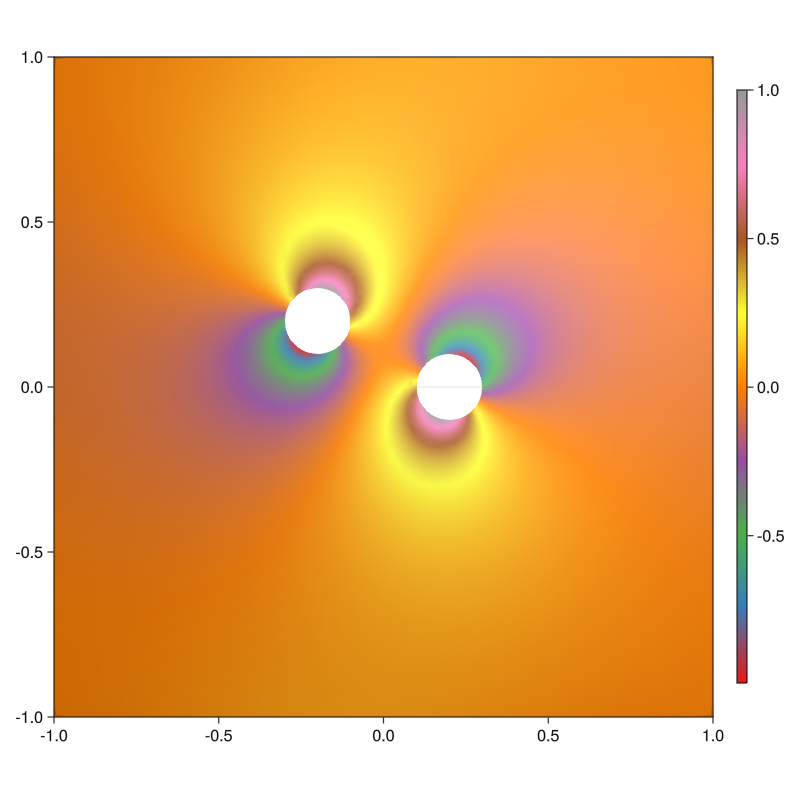}
  \includegraphics[width=0.28\textwidth]{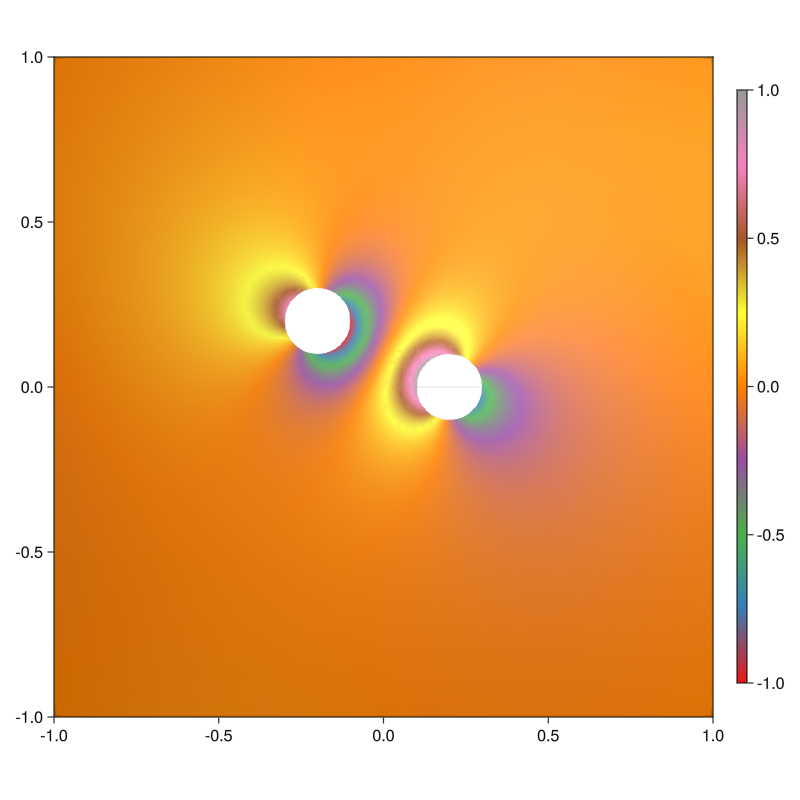}
  \includegraphics[width=0.28\textwidth]{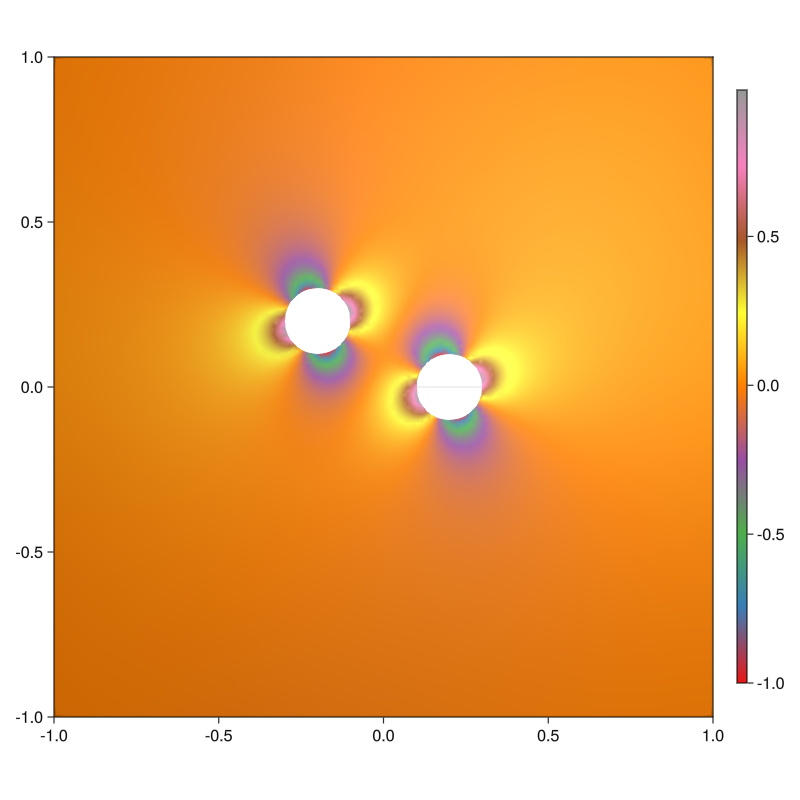}
  \caption{Approximate Steklov eigenfunctions  of indices $2$ to $7$  on a punctured square torus with two circular holes. See \cref{s:Steklov}.}
  \label{f:SteklovEigSquare2} 
\end{figure}

\begin{figure}
  \centering
  \includegraphics[width=0.28\textwidth]{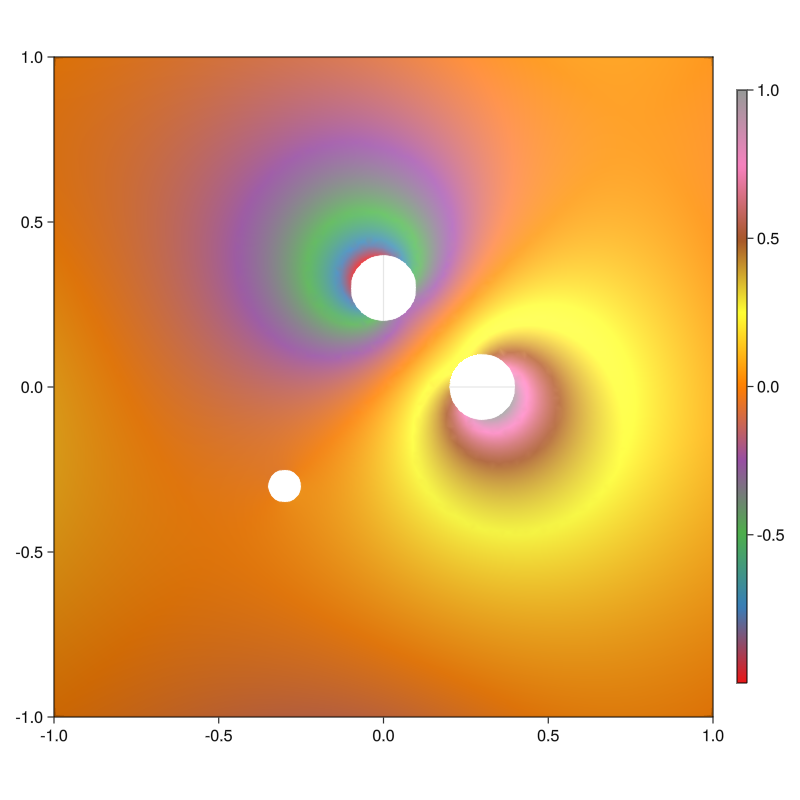}
  \includegraphics[width=0.28\textwidth]{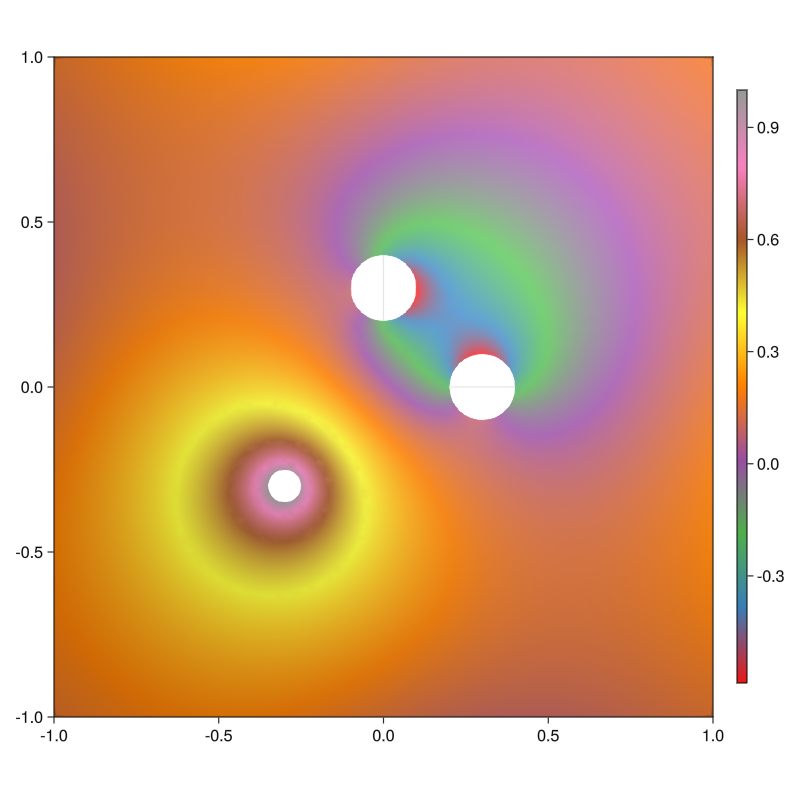}
  \includegraphics[width=0.28\textwidth]{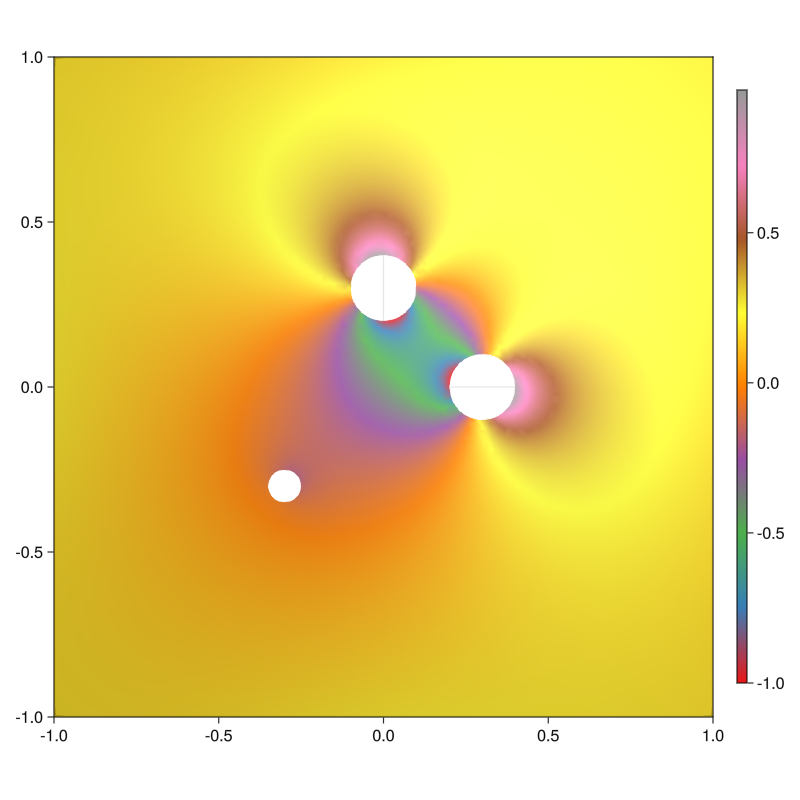}
  \includegraphics[width=0.28\textwidth]{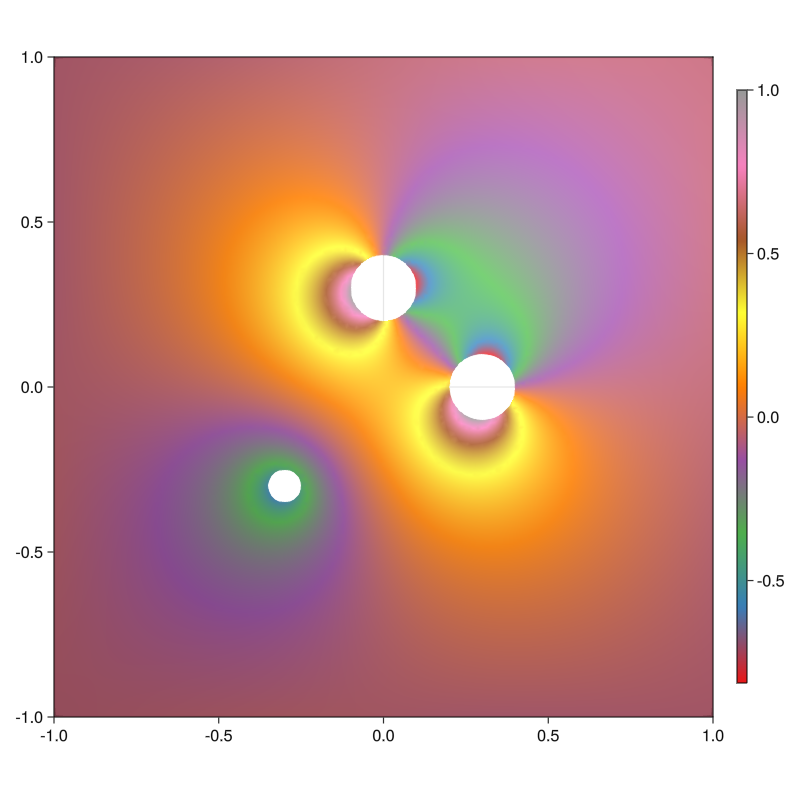}
  \includegraphics[width=0.28\textwidth]{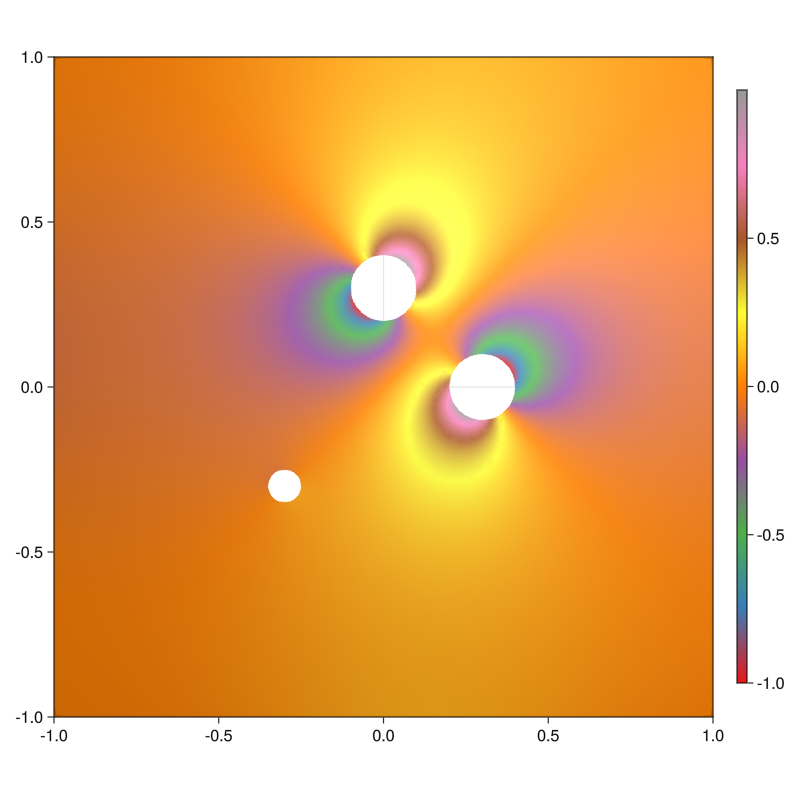}
  \includegraphics[width=0.28\textwidth]{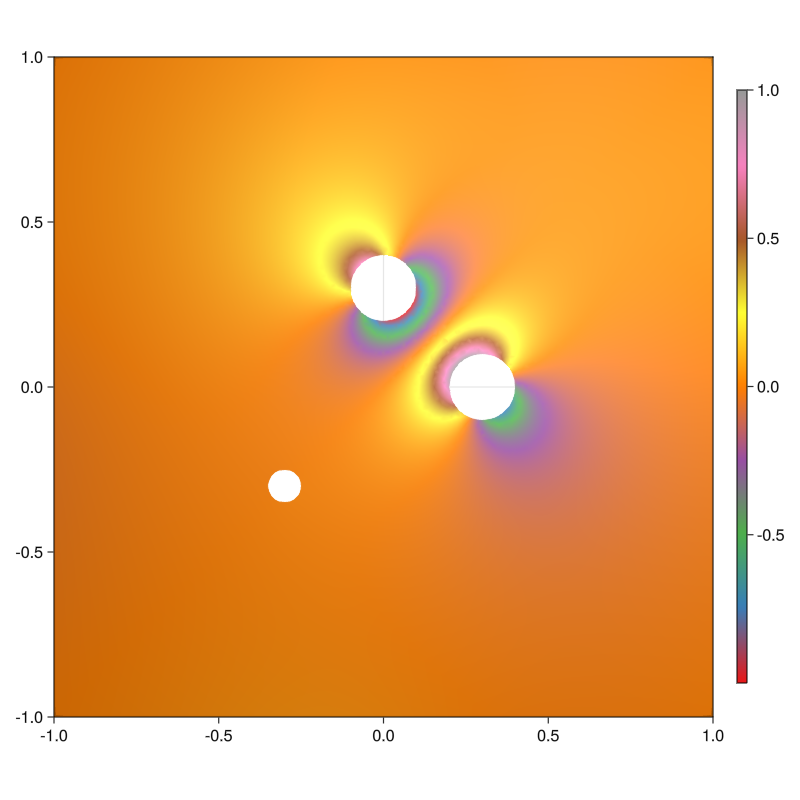}
  \caption{Approximate Steklov eigenfunctions  of indices $2$ to $7$  on a punctured square torus with three circular holes. See \cref{s:Steklov}.}
  \label{f:SteklovEigSquare3} 
\end{figure}

\begin{figure}[p]
  \centering
  \includegraphics[width=0.3\textwidth]{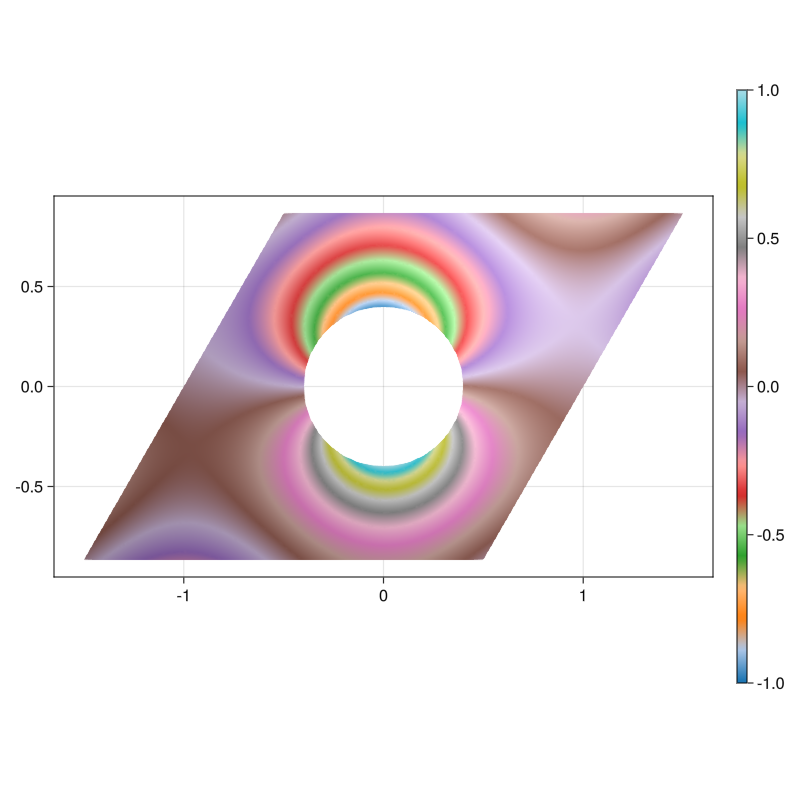}
  \includegraphics[width=0.3\textwidth]{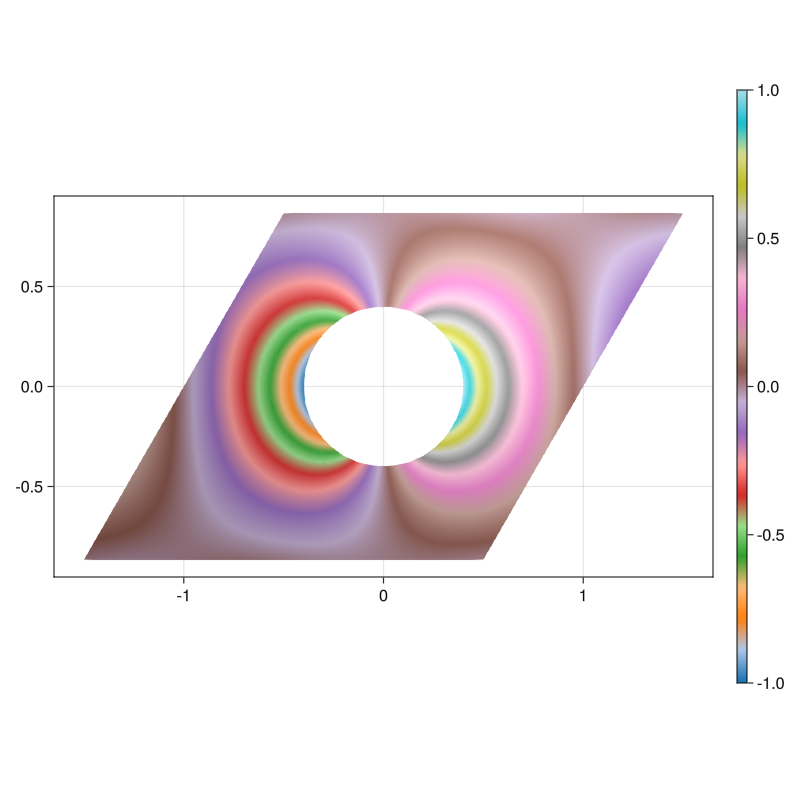}
  \includegraphics[width=0.3\textwidth]{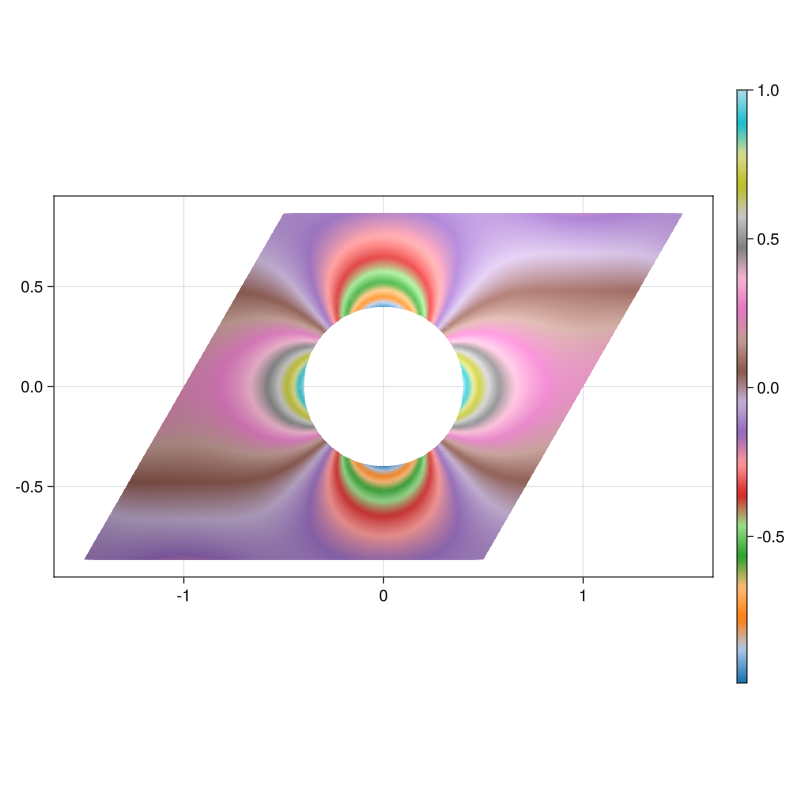}

  \vspace{-1.5cm}
  \includegraphics[width=0.3\textwidth]{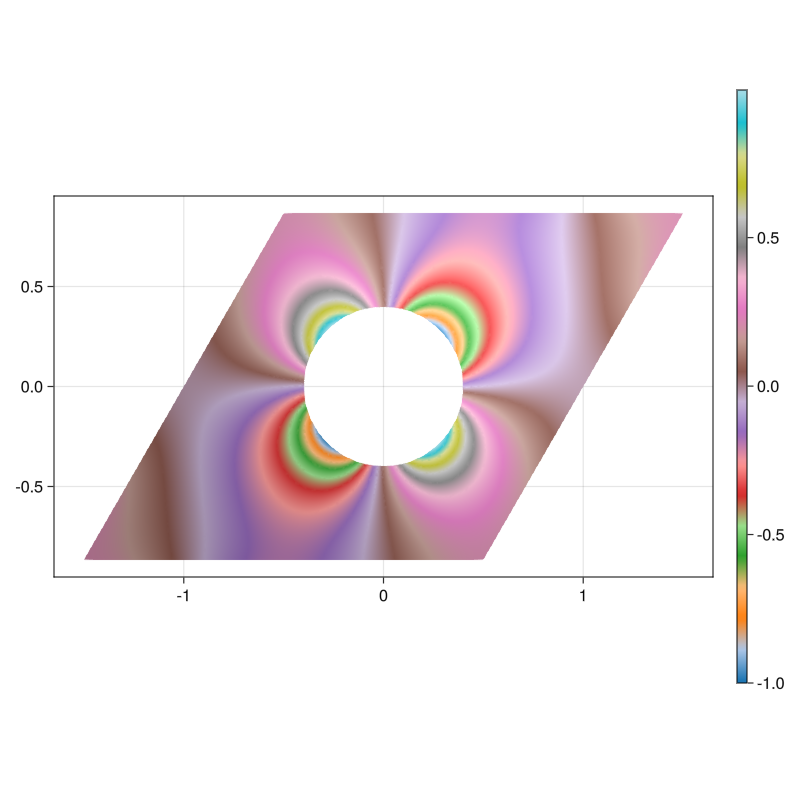}
  \includegraphics[width=0.3\textwidth]{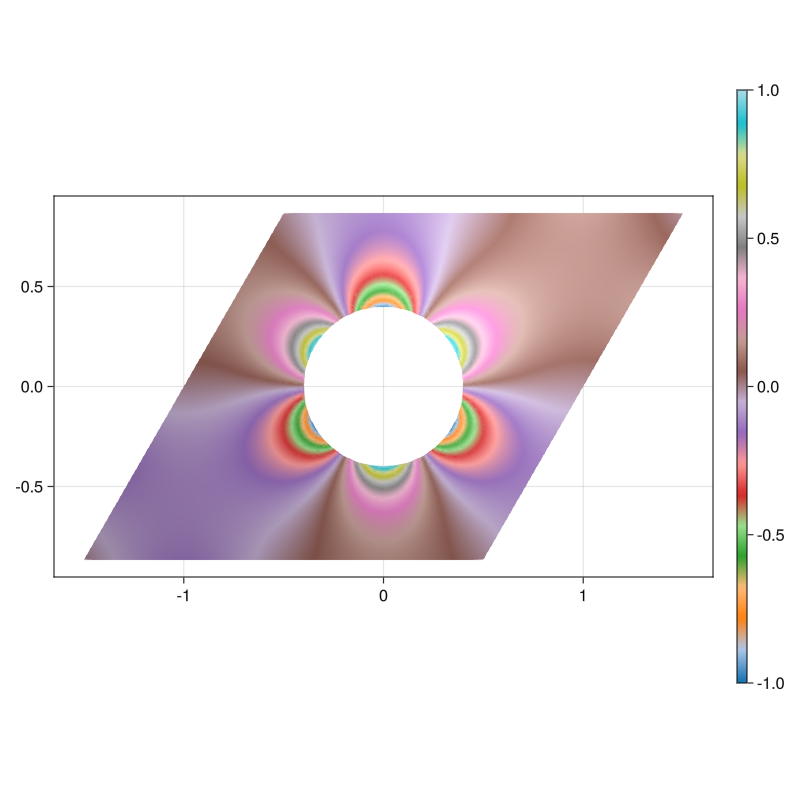}
  \includegraphics[width=0.3\textwidth]{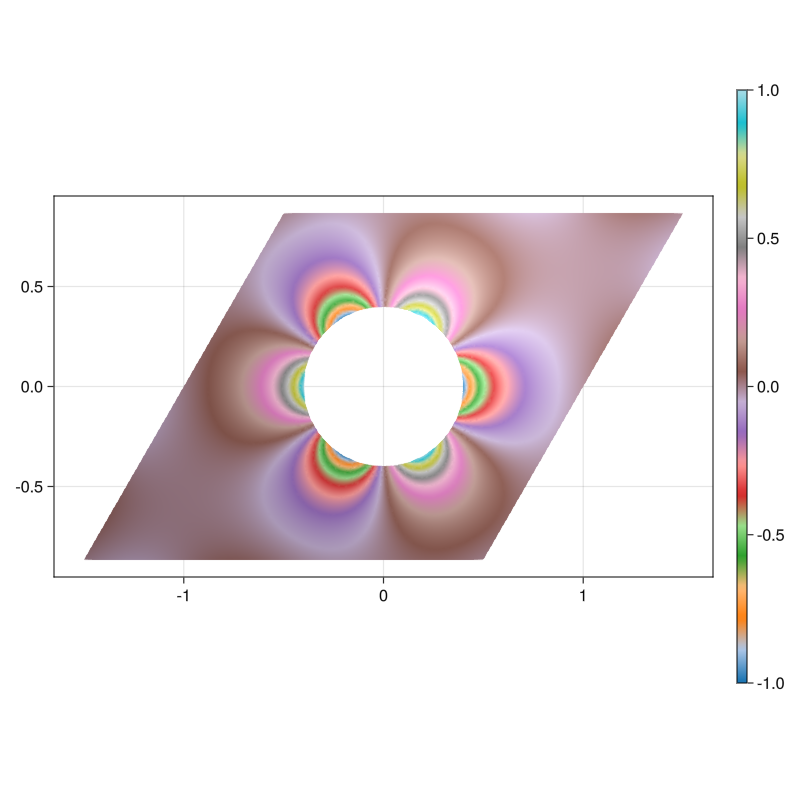}
  \caption{Approximate Steklov eigenfunctions of indices $2$ to $7$ on a punctured equilateral torus with one hole. See \cref{s:Steklov}.}
  \label{f:SteklovEigEqui1} 
\end{figure}

\begin{figure}
  \centering
  \includegraphics[width=0.3\textwidth]{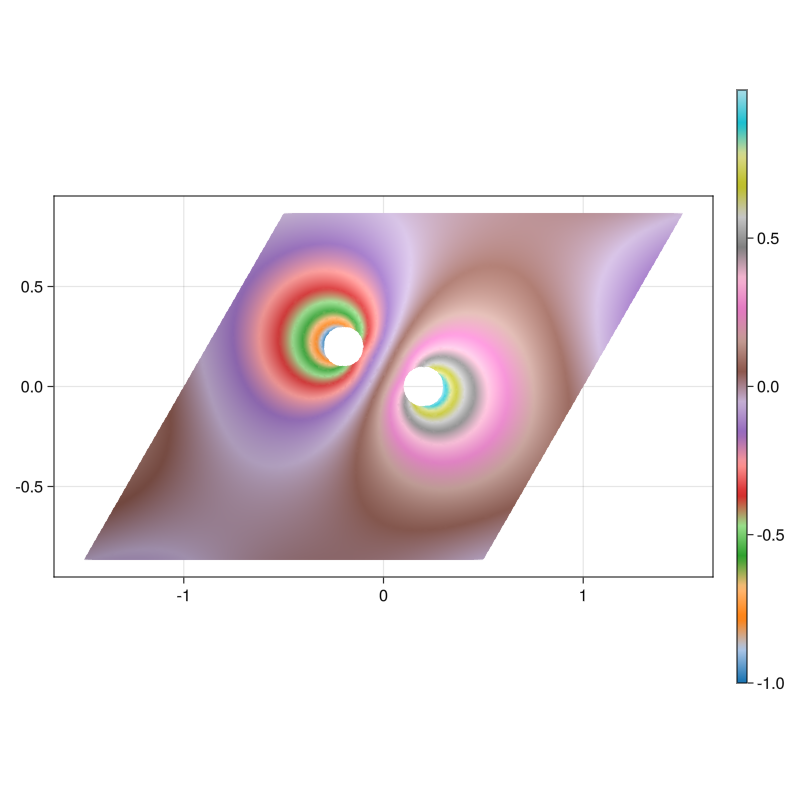}
  \includegraphics[width=0.3\textwidth]{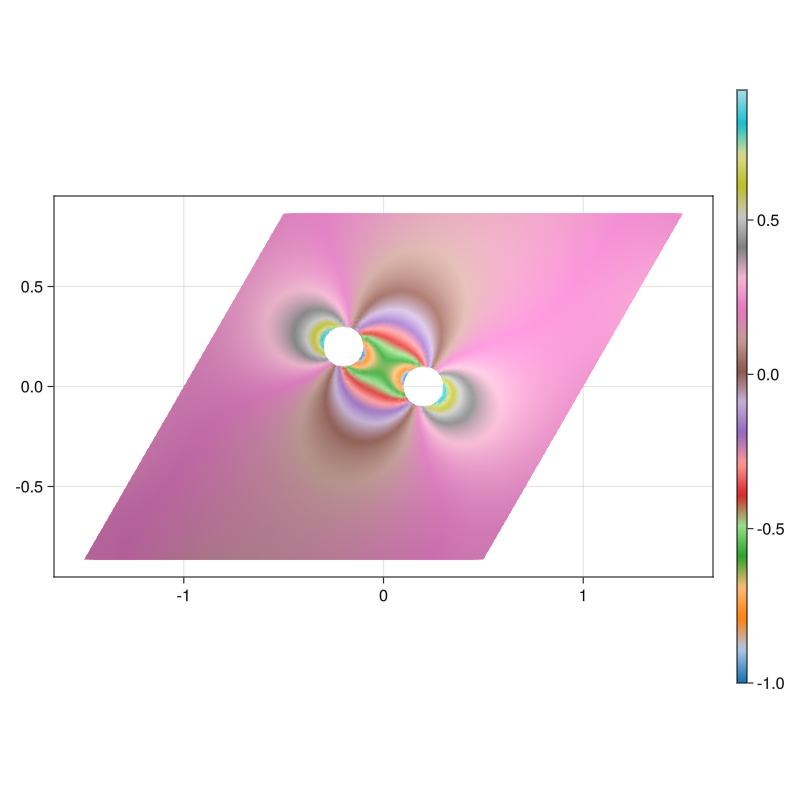}
  \includegraphics[width=0.3\textwidth]{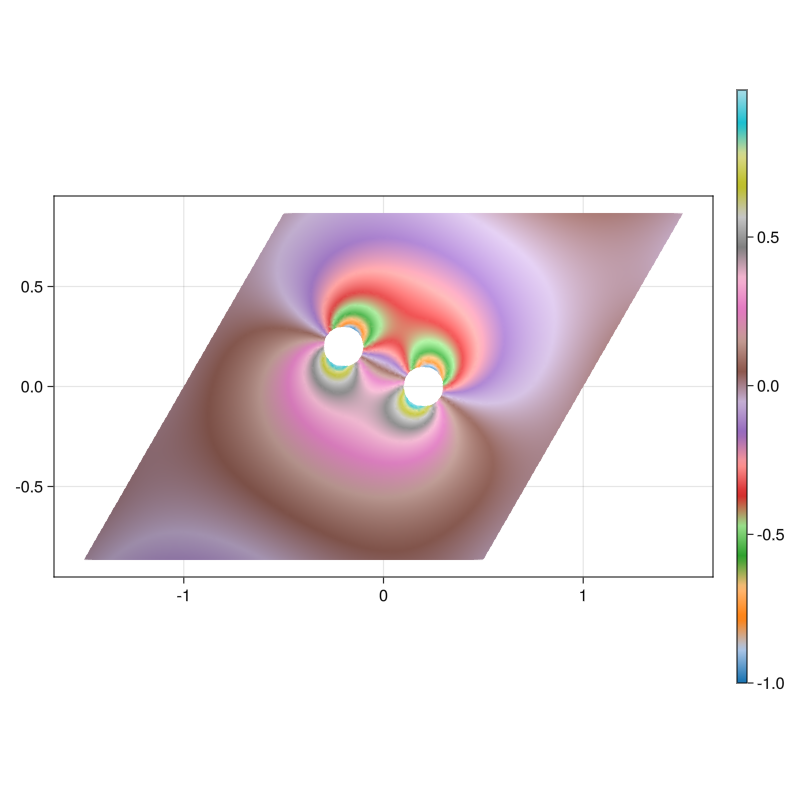}

  \vspace{-1.5cm}
  \includegraphics[width=0.3\textwidth]{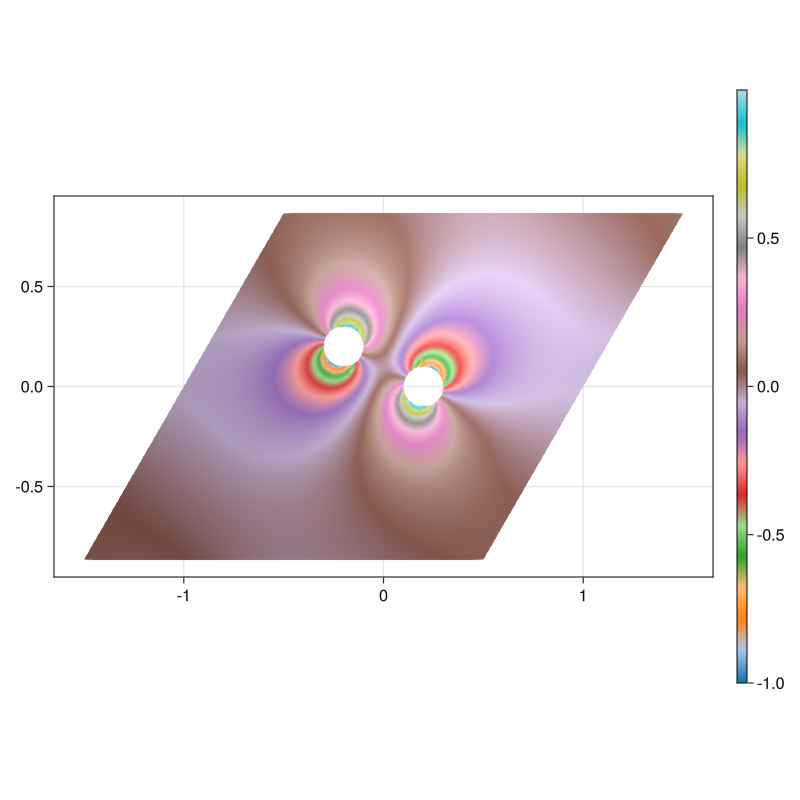}
  \includegraphics[width=0.3\textwidth]{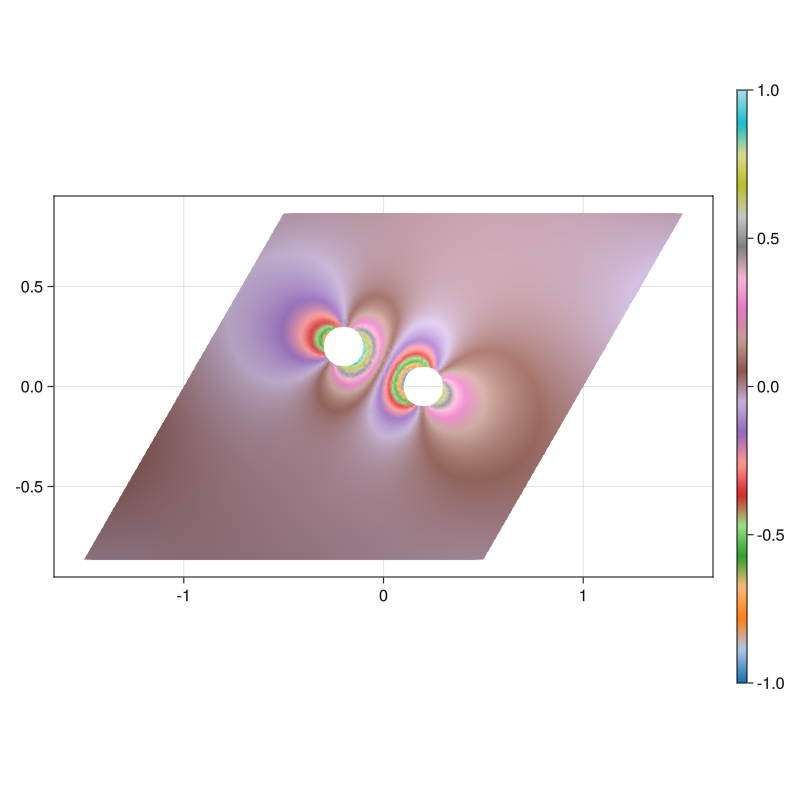}
  \includegraphics[width=0.3\textwidth]{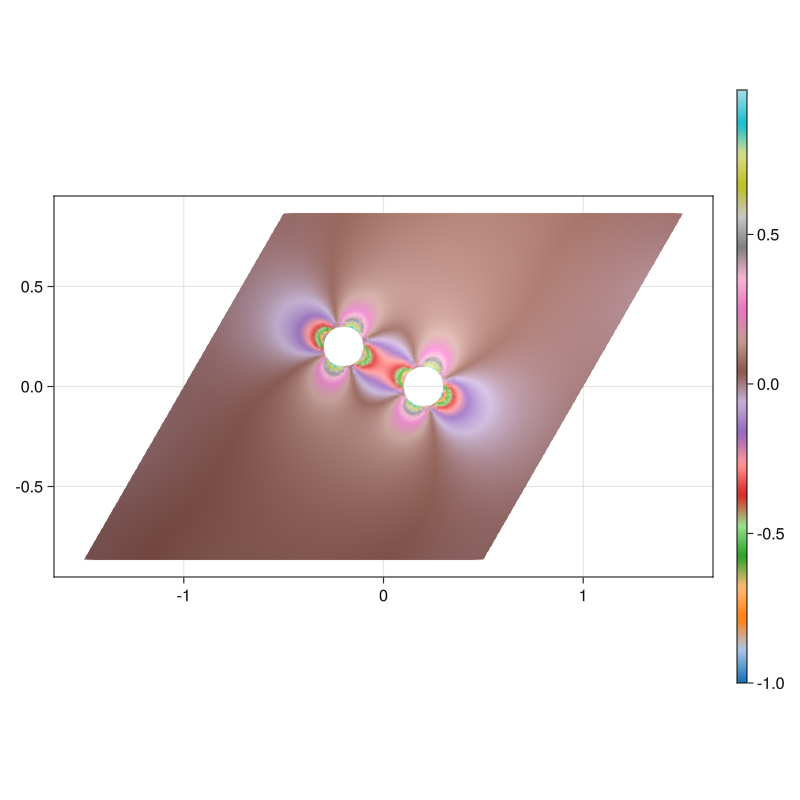}
  \caption{Approximate Steklov eigenfunctions  of indices $2$ to $7$  on a punctured equilateral torus with two circular holes. See \cref{s:Steklov}.}
  \label{f:SteklovEigEqui2} 
\end{figure}

\begin{figure}
  \centering
  \includegraphics[width=0.3\textwidth]{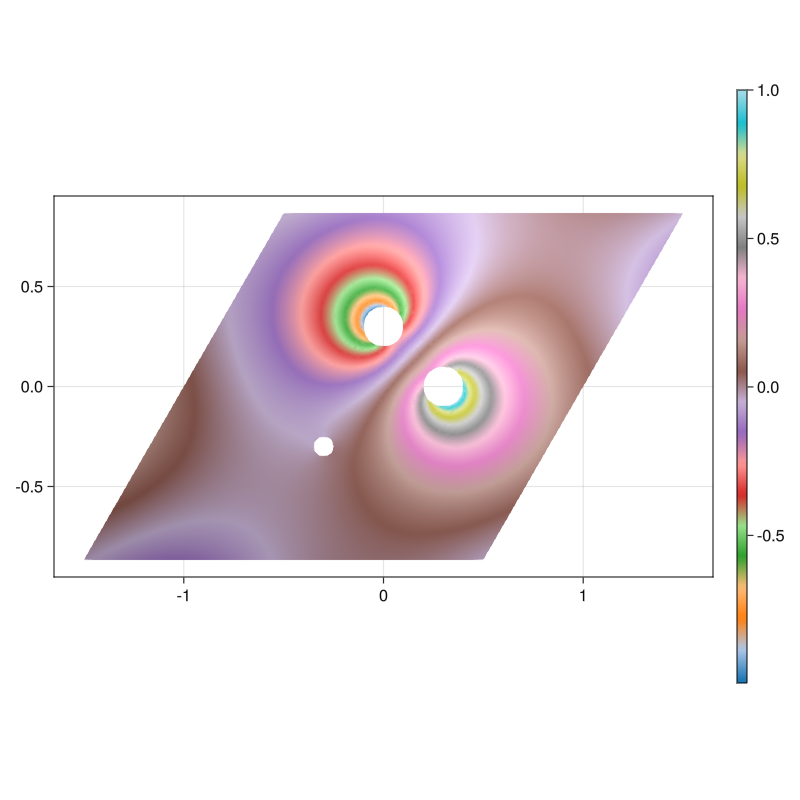}
  \includegraphics[width=0.3\textwidth]{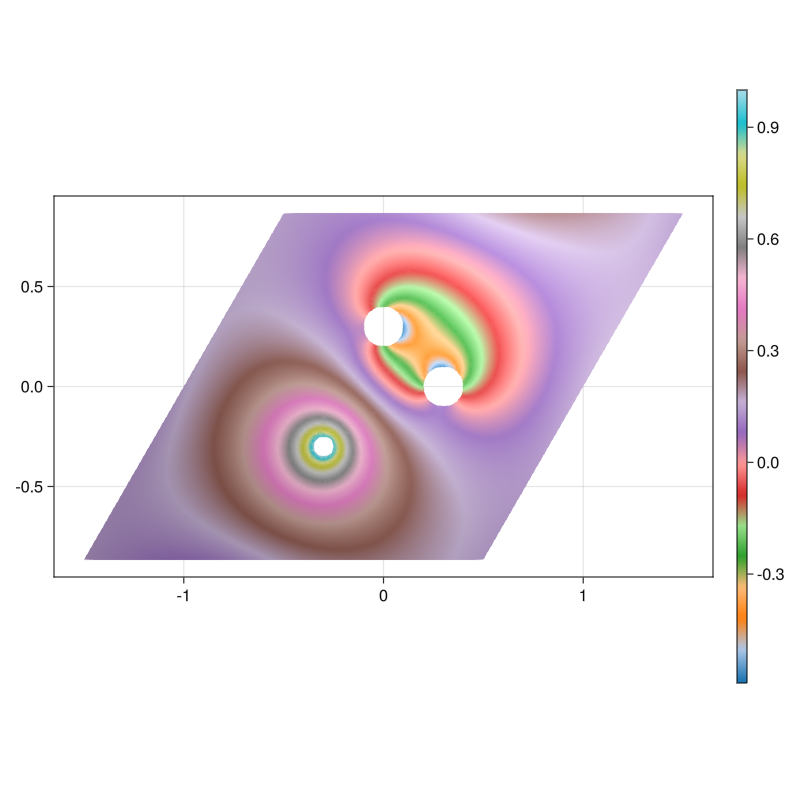}
  \includegraphics[width=0.3\textwidth]{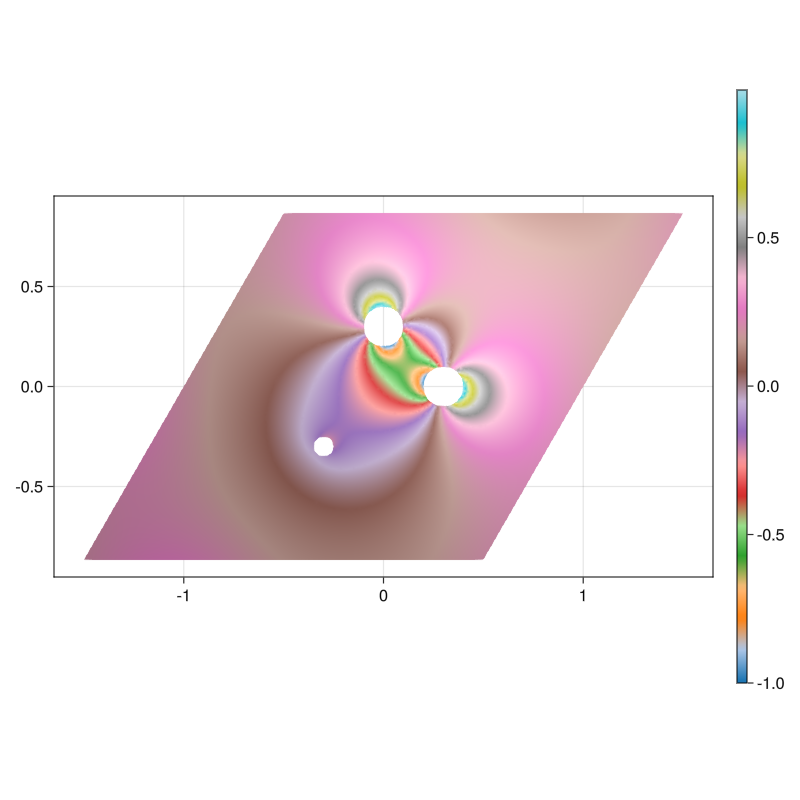}

  \vspace{-1.5cm}
  \includegraphics[width=0.3\textwidth]{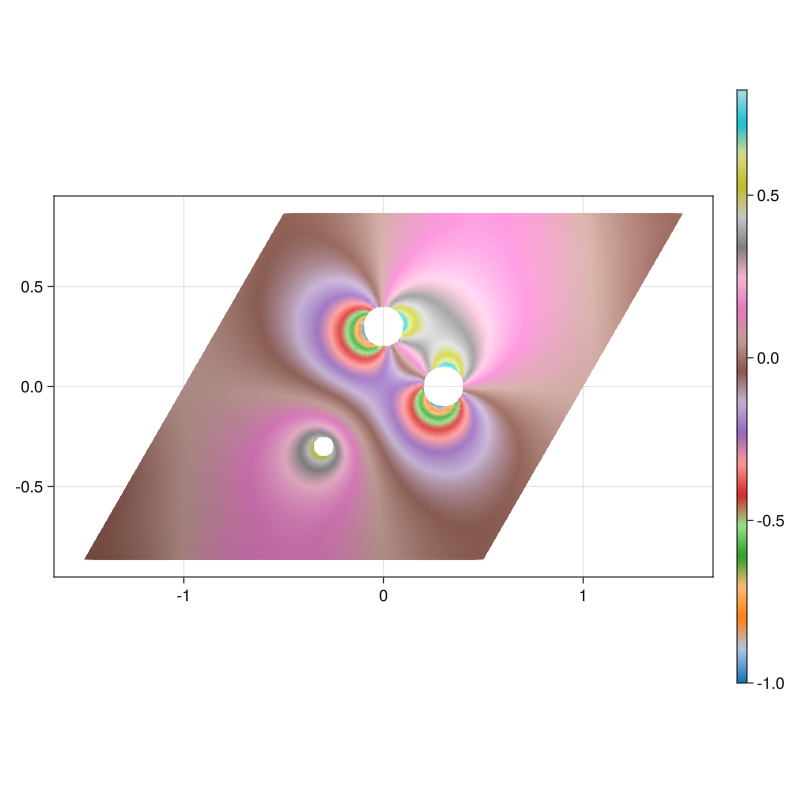}
  \includegraphics[width=0.3\textwidth]{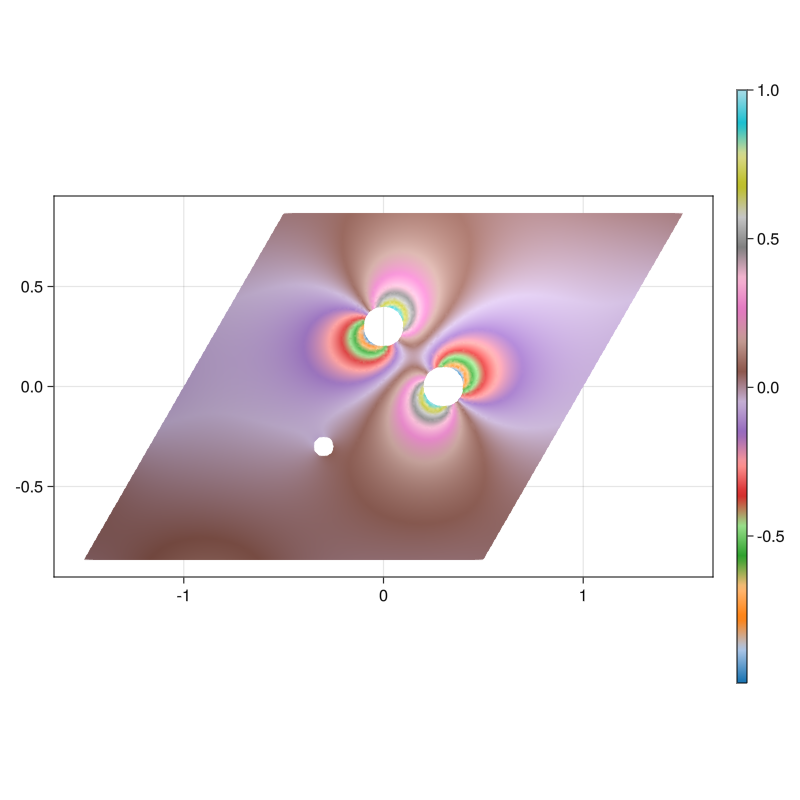}
  \includegraphics[width=0.3\textwidth]{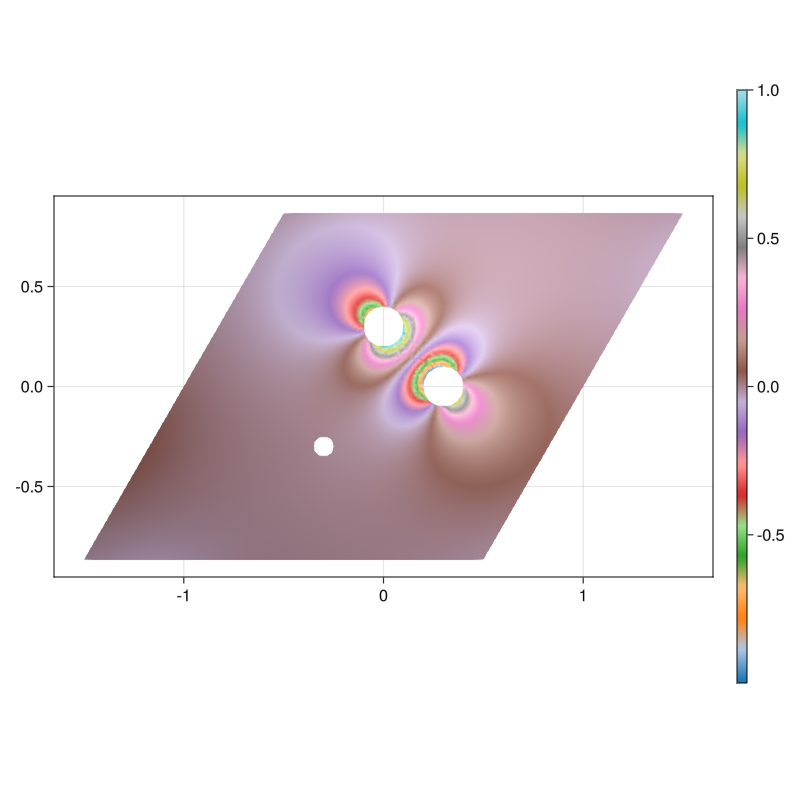}
  \caption{Approximate Steklov eigenfunctions  of indices $2$ to $7$  on a punctured equilateral torus with three circular holes. See \cref{s:Steklov}.}
  \label{f:SteklovEigEqui3} 
\end{figure}

\begin{figure}
\centering
\includegraphics[width=0.3\textwidth]{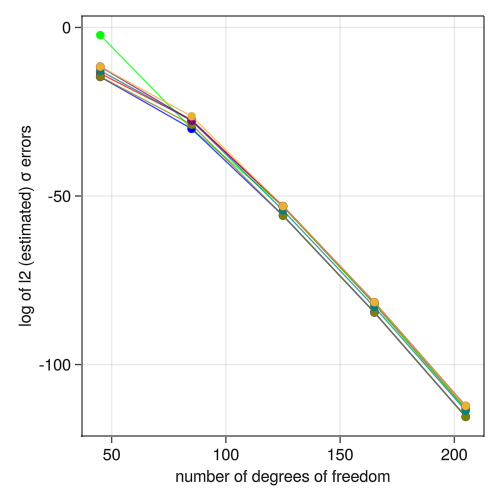}
\includegraphics[width=0.3\textwidth]{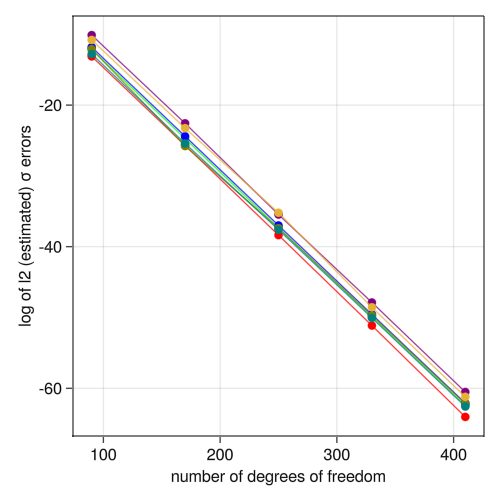}
\includegraphics[width=0.3\textwidth]{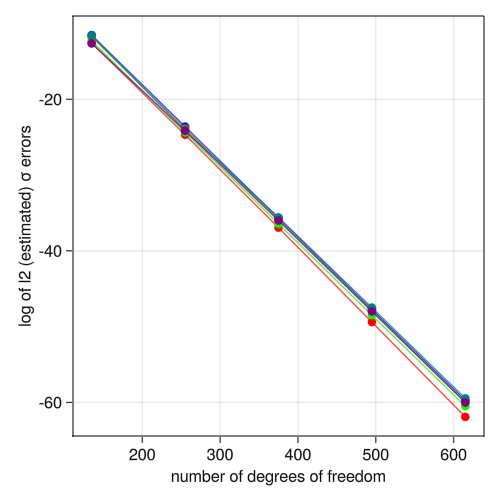}
\caption{Convergence plots for Steklov eigenvalues on a square domain with 1, 2 and 3 punctured holes. Each line corresponds to one of the first seven eigenvalues. See \cref{s:Steklov}.}
\label{f:SteklovConv}
\end{figure}

We study three geometrical configurations: tori which are the complement of 
$K_1 = B(a_1,r)$ with $a_1 = 0$ and $r=0.4$, the complement of $K_i = B(a_i,r_i)$, $i=1,2$ with $a_1 = 0.2$, $a_2 = -0.2+0.2i$ and $r_1 = r_2 = 0.1$ and the complement of $K_i = B(a_i,r_i)$, $i=1,2,3$ with $a_1 = 0.3$, $a_2 = 0.3i$, $a_3 = -0.3-0.3i$ and $r_1 = r_2 = 0.1$, $r_3 = 0.05$. We approximated Steklov  eigenfunctions of the square torus with half-periods $(\omega_1,\omega_2) = (1,i)$ (see   \cref{f:SteklovEigSquare1,f:SteklovEigSquare2,f:SteklovEigSquare3}) and of  the equilateral torus with half-periods $(\omega_1,\omega_2) = \left(1,\frac{1}{2} + \frac{\sqrt{3}}{2} \imath \right)$ in these three configurations (see  \cref{f:SteklovEigEqui1,f:SteklovEigEqui2,f:SteklovEigEqui3}).
The first eigenvalue is zero which corresponds to a constant eigenfunction. In these figures, Steklov eigenfunctions of indices 2 to 7 are plotted. The Steklov eigenfunctions, as expected, are oscillatory near the boundary and decay exponentially away from the boundary. We used \cref{p:steklovbound} to estimate the approximation error of the Steklov eigenvalues. We approximated the $L^2$ boundary term by a (periodic) trapezoidal quadrature formula after doubling the number of sampled points. Convergence plots for Steklov eigenvalues on a square domain with 1, 2, and 3 punctured circular holes are given in \cref{f:SteklovConv}. As expected, spectral convergence is also observed in these situations. The same convergence rate has also been obtained studying the equilateral case. 

Finally, in \cref{s:Eigenvalues}, we report in  \cref{t:Tab1,t:Tab2,t:Tab3,t:Tabequi1,t:Tabequi2,t:Tabequi3}  our approximation of the first six non-trivial eigenvalues of the square and equilateral tori with $b=1$, 2, and 3 circular holes. We believe that the reported $50$ digits are correct in each case.
As indicated in \cref{t:LogConjThmTorus}, when there is only one connected boundary component ($b=1$), the eigenfunctions do not involve the logarithmic term and  are oscillatory along the boundary as shown in  \cref{f:SteklovEigSquare1,f:SteklovEigEqui1}. In general, eigenfunctions corresponding to larger Steklov eigenvalues are more oscillatory near the boundary. Note that the tori parameters, $\omega$, effects the multiplicity of the eigenvalues. On a square torus with one circular hole, $\sigma_2=\sigma_3$ and $\sigma_6=\sigma_7$ while, on an equilateral torus with one circular hole, $\sigma_2=\sigma_3$ and $\sigma_4=\sigma_5$. Since the domains with two or three circular holes do not possess symmetry, we observe that all the obtained eigenvalues are simple.

\section{Discussion} \label{s:disc}
In this paper, we established \cref{t:LogConjThmTorus}, a Logarithmic Conjugation Theorem on finitely-connected tori. We used the theorem to find a series solution representation of harmonic functions on finitely-connected tori; see \cref{t:HarmonicPerDecomp}. Implementing the numerical method in Julia using arbitrary precision, we approximate solutions to the Laplace problem \eqref{e:Laplace} and the Steklov eigenvalue problem \eqref{e:Steklov}; see \cref{s:examples}. Using a posteriori estimation, we show that the approximate solution of the Laplace problem has error less than $10^{-100}$ using a few hundred degrees of freedom and the Steklov eigenvalues have similar error. 

There are several future directions for this work. 
The fundamental solution of Laplacian on flat tori can be expressed as a logarithmic function involving first Jacobi theta function \cite{Lin_2010,mamode2014fundamental}; we think it would be interesting to develop integral equation methods to approximate harmonic functions on finitely-connected tori in the spirit of \cite{barnett2018unified}.
We have focused on the case where the domain complement,  $ \cup_{j \in [b]} K_j$ has smooth boundary. We think it would be interesting to extend the methods in \cite{Gopal_2019} to improve the order of convergence for non-smooth boundaries. 
Finally, we think it would be interesting to apply the developed numerical methods to 
the numerical problem of computing extremal Steklov eigenvalue problems for finitely-connected flat tori \cite{Kao_2022,Kao_2020}.

\bibliographystyle{siam}
\bibliography{refs}

\clearpage
\appendix
\section{Computing normal derivatives} \label{s:NormalDerivatives}
In this appendix, we provide some details for computing normal derivatives of functions of a complex variable. 
Denote $f(z)=u(x,y)+\imath v(x,y)$ with $z = x + \imath y$. Since $f$ is analytic, we
have $u_{x}=v_{y}$ and $u_{y}=-v_{x}.$ Furthermore, $f_{x}=f'(z)$
and $f_{y}=\imath f'(z).$ Thus, with $n=n_{1}+\imath n_{2}$, we have 
\begin{align*}
u_{n} & = n_{1}u_{x}+n_{2}u_{y}
  =n_{1}u_{x}-n_{2}v_{x}
  =\Re \left[(n_{1}+\imath n_{2})(u_{x}+\imath v_{x})\right]
  =\Re (n f'(z)) \\ 
v_n & =n_{1}v_{x}+n_{2}v_{y}
 =-n_{1}u_{y}+n_{2}v_{y}
 =-\Im (n_{1}+\imath n_{2})(u_{y}+\imath v_{y})
 =\Im (n f'(z)). 
\end{align*}
For example, $f(z)=z^{k}$, 
\begin{align*}
\left(\Re \left(z^{k}\right)\right)_{n} & =k\Re (n z^{k-1}), \\
\left(\Im \left(z^{k}\right)\right)_{n} & =k\Im (n z^{k-1}).
\end{align*}
If $f(z)=\wp^{(k)}(z-a_{j}),$ 
\begin{align*}
\left(\Re \left(\wp^{(k)}(z-a_{j})\right)\right)_{n}
& =\Re (n \wp^{(k+1)}(z-a_{j})), \\
\left(\Im \left(\wp^{(k)}(z-a_{j})\right)\right)_{n}
& = \Im (n \wp^{(k+1)}(z-a_{j})).
\end{align*}
If $f(z)=\hat{\zeta}(z-a_{j}),$ 
\begin{align*}
\left(\Re \left(\hat{\zeta}\left(z-a_{j}\right)\right)\right)_{n} & =\left(\Re \left(\zeta\left(z-a_{j}\right)\right)\right)_{n}-\left(\Re \left(\gamma_{2}\left(z-a_{j}\right)+\frac{\pi}{A}\left(z-a_{j}\right)^{*}\right)\right)_{n}\\
 & =-\Re \left(n \wp\left(z-a_{j}\right)\right)-\left(\left(\Re \left(\gamma_{2}\right)+\frac{\pi}{A}\right)n_{1}-\Im \left(\gamma_{2}\right)n_{2}\right), \\ 
\left(\Im \left(\hat{\zeta}\left(z-a_{j}\right)\right)\right)_{n} & =\left(\Im \left(\zeta\left(z-a_{j}\right)\right)\right)_{n}-\left(\Im \left(\gamma_{2}\left(z-a_{j}\right)+\frac{\pi}{A}\left(z-a_{j}\right)^{*}\right)\right)_{n}\\
 & =-\Im \left(n \wp\left(z-a_{j}\right)\right)-\left(\Im \left(\gamma_{2}\right)n_{1}+\left(\Re \left(\gamma_{2}\right)-\frac{\pi}{A}\right)n_{2}\right),
\end{align*}
If $f(z)=\log\left|\hat{\sigma}\left(z-a_{j}\right)\right|,$ 
\begin{align*}
\left(\log\left|\hat{\sigma}\left(z-a_{j}\right)\right|\right)_{n} & =n_{1}\left(\log\left|\hat{\sigma}\left(z-a_{j}\right)\right|\right)_{x}+n_{2}\left(\log\left|\hat{\sigma}\left(z-a_{j}\right)\right|\right)_{y}\\
 & =n_{1}\left(\Re \left(-\frac{1}{2}\gamma_{2}z^{2}-\frac{1}{2}\frac{\pi}{A}|z|^{2}\right)\right)_{x}+n_{2}\left(\Re \left(-\frac{1}{2}\gamma_{2}z^{2}-\frac{1}{2}\frac{\pi}{A}|z|^{2}\right)\right)_{y}\\
 & \quad +n_{1}\left(\log\left|\sigma\left(z-a_{j}\right)\right|\right)_{x}+n_{2}\left(\log\left|\sigma\left(z-a_{j}\right)\right|\right)_{y}\\
 & =n_{1}\left(\left(-\Re \left(\gamma_{2}\right)-\frac{\pi}{A}\right)x+\Im \left(\gamma_{2}\right)y\right)+n_{2}\left(\Im \left(\gamma_{2}\right)x+\left(\Re \left(\gamma_{2}\right)-\frac{\pi}{A}\right)y\right)\\
 & \quad +n_{1}\left(\Re \left(\zeta\left(z-a_{j}\right)\right)\right)+n_{2}\left(-\Im \left(\zeta\left(z-a_{j}\right)\right)\right)
\end{align*}

\clearpage
\section{Numerical values of computed Steklov eigenvalues}
\label{s:Eigenvalues}

\noindent Values of computed Steklov eigenvalues are given in \cref{t:Tab1,t:Tab2,t:Tab3,t:Tabequi1,t:Tabequi2,t:Tabequi3}; see \cref{s:Steklov} for details. 

\vspace{-.2cm}

{\small

\begin{center}
  \begin{table}[H]
    \begin{tabular}{|c|r|} 
      \hline
      $\sigma_2$  & 3.21737540790552735473880286001400036767774798208487 \\ 
      $\sigma_3$  & 3.21737540790552735473880286001400036767774798208487 \\ 
      $\sigma_4$  & 4.85099530552467697892257589130439715581461931719259 \\
      $\sigma_5$  & 5.15358084940676223549771471754234765157435969419525 \\
      $\sigma_6$  & 7.50305008416767542642635086056165243882709526430554 \\
      $\sigma_7$  & 7.50305008416767542642635086056165243882709526430554 \\
      \hline
    \end{tabular}
    \vspace{-0.3cm}
    \caption{Steklov eigenvalues of a square torus with one circular hole.}
    \label{t:Tab1} 
  \end{table}

\vspace{-0.3cm}

  \begin{table}[H]
    \begin{tabular}{ |c|r|} 
      \hline
      $\sigma_2$  & 6.45837308842285506198400983365912091999317179119988   \\
      $\sigma_3$  & 9.04038374077713587651429965380130970292955686420981   \\
      $\sigma_4$  & 9.32931391918711635803886895114746515357566095450257   \\
      $\sigma_5$  & 11.02561512617948586321622981756262835104523220458063  \\
      $\sigma_6$  & 12.69568331719729045908212485186369130848598658103989  \\
      $\sigma_7$  & 19.72884655790348748027382339516572459572547368362522  \\
      \hline
    \end{tabular}
    \vspace{-0.3cm}
    \caption{Steklov eigenvalues of a square torus with two circular holes.}
    \label{t:Tab2} 
  \end{table}

\vspace{-0.3cm}

  \begin{table}[H]
    \begin{tabular}{ |c|r|} 
      \hline
      $\sigma_2$  & 6.54721983775026738598476089606442586801693676638247  \\
      $\sigma_3$  & 6.79298688602543949226783518103096724408533776232952  \\ 
      $\sigma_4$  & 9.02715360305747386008778464587475727275979230551042  \\
      $\sigma_5$  & 9.75911376018587533254687022367601130464658416864329  \\
      $\sigma_6$  & 11.11563661826511047191549742109769883301063010901993 \\
      $\sigma_7$  & 13.08067309361125105561475152956318658177620553096385 \\
      \hline
    \end{tabular}
        \vspace{-0.3cm}
    \caption{Steklov eigenvalues of a square torus with three circular holes}
    \label{t:Tab3} 
  \end{table}

\vspace{-0.3cm}

  \begin{table}[H]
    \begin{tabular}{ |c|r|} 
      \hline
      $\sigma_2$ & 3.34865594380260534169550288243470971962587318064277  \\
      $\sigma_3$ & 3.34865594380260534169550288243470971962587318064277  \\
      $\sigma_4$ & 4.99978881548382813234141616969113198885117552416465  \\
      $\sigma_5$ & 4.99978881548382813234141616969113198885117552416465  \\
      $\sigma_6$ & 7.44392530690947308002824485738760008901145380307620  \\
      $\sigma_7$ & 7.55649710043624518482844840631875099119732734059433  \\
      \hline
    \end{tabular}
    \vspace{-0.3cm}
    \caption{Steklov eigenvalues of an equilateral torus with one hole.}
    \label{t:Tabequi1} 
  \end{table}

\vspace{-0.3cm}

  \begin{table}[H]
    \begin{tabular}{ |c|r|} 
      \hline
      $\sigma_2$  & 6.53794803818597918794030349125758145344842633243163  \\
      $\sigma_3$  & 9.03760803330365503342990995931942991592541389841134  \\
      $\sigma_4$  & 9.37148419781059159007737134528684902568383756667729  \\
      $\sigma_5$  & 11.02904931936017784776119216982004095594847249520813 \\
      $\sigma_6$  & 12.70222698966325001285792418382443547595163064198654 \\
      $\sigma_7$  & 19.72940718569248148461882657324755544321541234433839 \\
      \hline
    \end{tabular}
        \vspace{-0.3cm}
    \caption{Steklov eigenvalues of an equilateral torus with two circular holes.}
    \label{t:Tabequi2} 
  \end{table}

\vspace{-0.3cm}

  \begin{table}[H]
    \begin{tabular}{ |c|r|} 
      \hline
      $\sigma_2$  & 6.63530737085667505246439432756580077469480498850424  \\
      $\sigma_3$  & 6.94424494471680808970061612948991192950478141474806  \\
      $\sigma_4$  & 9.02318420302479178183722837786227147321263458092783  \\
      $\sigma_5$  & 9.69311259795549433304394048074564041975314036318590  \\
      $\sigma_6$  & 11.14183481942696624006786357768349746325365435641531 \\
      $\sigma_7$  & 13.09270988086125229485281063758500867332548835687565 \\
      \hline
    \end{tabular}
        \vspace{-0.3cm}
    \caption{Steklov eigenvalues of an equilateral torus with three circular holes.}     
    \label{t:Tabequi3}
  \end{table}

\end{center}}

\end{document}